\documentclass[11pt,leqno]{amsart}
\usepackage{epsfig}
\usepackage{amssymb}
\usepackage{amscd}
\usepackage[matrix,arrow]{xy}

\newtheorem{theo}{Theorem}[section]
\newtheorem{lemma}[theo]{Lemma}
\newtheorem{defi}[theo]{Definition}
\newtheorem{prop}[theo]{Proposition}

\newtheorem{cor}[theo]{Corollary}
\newtheorem{remark}[theo]{Remark}

\newtheorem{example}[theo]{Example}
\numberwithin{equation}{section}

\def\B{{\mathcal B}}

\def\P{{\mathcal{P}}}
\def\D{{\mathcal{D}}}

\def\bR{{\mathbf R}}
\def\bL{{\mathbf L}}

\def\pre-tr{\operatorname{pre-tr}}
\def\h{\operatorname{h}}
\def\Hom{\operatorname{Hom}}
\def\End{\operatorname{End}}
\def\gr{\operatorname{gr}}

\newcommand{\bbZ}{{\mathbb Z}}

\newcommand{\cJ}{{\mathcal J}}
\newcommand{\cQ}{{\mathcal Q}}

\newcommand{\cO}{{\mathcal O}}
\newcommand{\cP}{{\mathcal P}}

\newcommand{\cM}{{\mathcal M}}

\newcommand{\cD}{{\mathcal D}}

\newcommand{\cA}{{\mathcal A}}
\newcommand{\cB}{{\mathcal B}}
\newcommand{\cI}{{\mathcal I}}
\newcommand{\cC}{{\mathcal C}}
\newcommand{\cE}{{\mathcal E}}

\newcommand{\cR}{{\mathcal R}}

\newcommand{\cl}{\operatorname{cl}}
\newcommand{\qu}{\operatorname{qu}}

\newcommand{\DG}{\operatorname{DG}}
\newcommand{\Fun}{\operatorname{Fun}}

\newcommand{\Def}{\operatorname{Def}}

\newcommand{\Ker}{\operatorname{Ker}}
\newcommand{\im}{\operatorname{Im}}

\newcommand{\Ext}{\operatorname{Ext}}

\newcommand{\Id}{\operatorname{Id}}

\newcommand{\Ind}{\operatorname{Ind}}
\newcommand{\Res}{\operatorname{Res}}

\newcommand{\dgart}{\operatorname{dgart}}
\newcommand{\art}{\operatorname{art}}

\newcommand{\coDef}{\operatorname{coDef}}
\newcommand{\cart}{\operatorname{cart}}
\newcommand{\Alg}{\operatorname{Alg}}

\newcommand{\Ho}{\operatorname{Ho}}

\newcommand{\id}{\operatorname{id}}

\newcommand{\dgalg}{\operatorname{dgalg}}
\newcommand{\DEF}{\operatorname{DEF}}
\newcommand{\coDEF}{\operatorname{coDEF}}

\newcommand{\adgalg}{\operatorname{adgalg}}

\title[Deformation theory of objects  in homotopy and derived categories I]
{Deformation theory of objects  in homotopy and derived categories I: general theory}
\author{Alexander I.~Efimov}
\address{Department of Mechanics and Mathematics, Moscow State University, Moscow,
Russia and Independent Univeroty of Moscow}
\email{efimov@mccme.ru}
\author{Valery A.~Lunts}
\address{Department of Mathematics, Indiana University,
Bloomington, IN 47405, USA} \email{vlunts@indiana.edu}
\author{Dmitri O.~Orlov}
\address{Algebra Section, Steklov Mathematical Institute, 8 Gubkina str., Moscow, 119991 Russia}
\email{orlov@mi.ras.ru}

\thanks{The first
named author was partially supported by grant NSh-1983.2008.1. The
second named author was partially supported by the NSA grant
H98230-05-1-0050 and CRDF grant RUM1-2661-MO-05. The third named
author was partially supported by CRDF grant RUM1-2661-MO-05, grant
RFFI 05-01-01034, grant INTAS 05-1000008-8118 and grant
NSh-9969.2006.1}

\begin{document}

\begin{abstract} This is the first paper in a series. We develop a general deformation theory of objects
in homotopy and derived categories of DG categories. Namely, for a
DG module $E$ over a DG category we define four deformation functors
$\Def ^{\h}(E)$, $\coDef ^{\h}(E)$, $\Def (E)$, $\coDef (E)$. The
first two functors describe the deformations (and co-deformations)
of $E$ in the homotopy category, and the last two - in the derived
category. We study their properties and relations. These functors
are defined on the category of artinian (not necessarily
commutative) DG algebras.
\end{abstract}

\maketitle

\tableofcontents

\section{Introduction}

It is well known (see for example \cite{De1}, \cite{De2},
\cite{Dr2}, \cite{Ge1}, \cite{Ge2}, \cite{Hi}) that for many
mathematical objects $X$ (defined over a field of characteristic
zero) the formal deformation theory of $X$ is controlled by a DG Lie
algebra $\mathfrak{g}=\mathfrak{g}(X)$ of (derived) infinitesimal
automorphisms of $X$. This is so in case $X$ is an algebra, a
compact complex manifold, a principal $G$-bundle, etc..

Let $\cM(X)$ denote the base of the universal deformation of $X$ and
$o\in \cM(X)$ be the point corresponding to $X$. Then (under some
conditions on $\mathfrak{g}$) the completion of the local ring
$\hat{\cO}_{\cM(X),o}$ is naturally isomorphic to the linear dual of
the homology space $H_0(\mathfrak{g})$. The space
$H_0(\mathfrak{g})$ is a co-commutative coalgebra, hence its dual is
a commutative algebra.

The homology $H_0(\mathfrak{g})$ is the zero cohomology group of
$B\mathfrak{g}$ -- the bar construction of $\mathfrak{g}$, which is
a co-commutative DG coalgebra. It is therefore natural to consider
the DG "formal moduli space" $\cM ^{DG}(X)$, so that the
corresponding completion $\hat{\cO}_{\cM^{DG}(X),o}$ of the "local
ring"  is the linear dual $(B\mathfrak{g})^*$, which is a
commutative DG algebra. The space $\cM ^{DG}(X)$ is thus the "true"
universal deformation space of $X$; it coincides with $\cM(X)$ in
case $H^i(B\mathfrak{g})=0$ for $i\neq 0$.  In particular, it
appears that the primary object is not the DG algebra
$(B\mathfrak{g})^*$, but rather the DG coalgebra $B\mathfrak{g}$
(this is the point of view in \cite{Hi}). In any case, the
corresponding deformation functor is naturally defined on the
category of commutative artinian DG algebras (see \cite{Hi}).

Note that the passage from a DG Lie algebra $\mathfrak{g}$ to the
commutative DG algebra $(B\mathfrak{g})^*$ is an example of the
Koszul duality for operads \cite{GiKa}. Indeed, the operad of DG Lie
algebras is Koszul dual to that of commutative DG algebras.

Some examples of DG algebraic geometry are discussed in  \cite{Ka},
\cite{Ci-FoKa1}, \cite{Ci-FoKa2}.

This paper (and the following papers \cite{LOII}, \cite{LOIII}) is
concerned with a general deformation theory in a slightly different
context. Namely, we consider deformations of "linear" objects $E$,
such as objects in a homotopy or a derived category. More precisely,
$E$ is a right DG module over a DG category $\cA$. In this case the
deformation theory of $E$ is controlled by $\cB =\End (E)$ which is
a DG {\it algebra} (and not a DG Lie algebra). (This works equally
well in positive characteristic.) Then the DG formal deformation
space of $E$ is the "Spec" of the (noncommutative!) DG algebra
$(B\cB )^*$ -- the linear dual of the bar construction $B\cB$ which
is a DG coalgebra. Again this is in agreement with the Koszul
duality for operads, since the operad of DG algebras is self-dual.
(All this was already anticipated in \cite{Dr2}.)

More precisely, let $\dgart$ be the category of local artinian (not
necessarily commutative) DG algebras and $\bf{Gpd}$ be the
2-category of groupoids. For a right DG module $E$ over a DG
category $\cA $ we define four pseudo-functors
$$\Def ^{\h}(E), \coDef ^{\h}(E), \Def (E), \coDef (E):\dgart \to
{\bf Gpd}.$$ The first two are the {\it homotopy} deformation and
co-deformation pseudo-functors, i.e. they describe deformations (and
co-deformations) of $E$ in the homotopy category of DG $\cA
^{op}$-modules; and the last two are their {\it derived} analogues. We
prove that the pseudo-functors $\Def ^{\h}(E)$, $\coDef ^{\h}(E)$
are equivalent and depend only on the quasi-isomorphism class of the
DG algebra $\End (E)$. The derived pseudo-functors $\Def (E)$,
$\coDef (E)$ need some boundedness conditions to give the "right"
answer and in that case they are equivalent to $\Def ^{\h}(F)$ and $
\coDef ^{\h}(F)$ respectively for an appropriately chosen
h-projective or h-injective DG module $F$ which is quasi-isomorphic
to $E$ (one also needs to restrict the pseudo-functors to the
category $\dgart _-$ of negative artinian DG algebras).

This first paper is devoted to the study of general properties of
the above four pseudo-functors and relations between them. Part 1 of
the paper is a rather lengthy review of basics of DG categories and
DG modules over them with some minor additions that we did not find
in the literature. The reader who is familiar with basic DG
categories is suggested to go directly to Part 2, except for looking
up the definition of the DG functors $i^*$ and $i^!$.

In the second paper \cite{LOII} we study the pro-representability of
these pseudo-functors. Recall that "classically" one defines
representability only for functors with values in the category of
sets (since the collection of morphisms between two objects in a
category is a set). For example, given a moduli problem in the form
of a pseudo-functor with values in the 2-category of groupoids one
then composes it with the functor $\pi _0$ to get a set valued
functor, which one then tries to (pro-) represent. This is certainly
a loss of information. But in order to represent the original
pseudo-functor one needs the source category to be a bicategory.

It turns out that there is a natural bicategory $2\text{-}\adgalg$
of augmented DG algebras. (Actually we consider two versions of this
bicategory, $2\text{-}\adgalg$ and $2^\prime\text{-}\adgalg$, but
then show that they are equivalent). We consider its full
subcategory $2\text{-}\dgart _-$ whose objects are negative artinian
DG algebras, and show that the derived deformation functors can be
naturally extended to pseudo-functors
$$\coDEF _-(E):2\text{-}\dgart _- \to {\bf Gpd},\quad \DEF _-(E):2^\prime\text{-}\dgart _- \to {\bf
Gpd}.$$ Then (under some finiteness conditions on the cohomology
algebra $H(\cC)$ of the DG algebra $\cC=\bR \Hom (E,E)$) we prove
pro-representability of these pseudo-functors by some local
complete DG algebra described by means of $A_{\infty}$-structure
on $H(\cC)$.

This pro-representability appears to be more "natural" for the
pseudo-functor $\coDEF _-$, because there exists a "universal
co-deformation" of the DG $\cC^{op}$-module $\cC$. The
pro-representability of the pseudo-functor $\DEF _-$ may then be
formally deduced from that of $\coDEF _-$.

In the third paper \cite{LOIII} we show how to apply our deformation
theory of DG modules to deformations of complexes over abelian
categories. We also discuss examples from algebraic geometry.

We note that the noncommutative deformations (i.e. over
noncommutative artinian rings) of modules were already considered by
Laudal in \cite{La}. The basic difference between our work and
\cite{La} (besides the fact that our noncommutative artinian
algebras are DG algebras) is that we work in the derived context.
That is we only deform the differential in a suitably chosen complex
and keep the module structure constant.

It is our pleasure to thank A.Bondal, P.Deligne, M.Mandell, M.Larsen
and P.Bressler for useful discussions. We especially appreciate the
generous help of B.Keller. We also thank W.Goldman and V.Schechtman
for sending us copies of letters \cite{De1} and \cite{Dr2}
respectively and W.Lowen for sending us the preprint \cite{Lo}. We
also thank J.Stasheff for his useful comments on the first version
of this paper.

\part{Preliminaries on DG categories}

\section{Artinian DG algebras}

 We fix a field $k$. All algebras are assumed to be
$\bbZ$ graded $k$-algebras with unit and all categories are
$k$-linear. Unless mentioned otherwise $\otimes $ means $\otimes
_k$.

For a homogeneous element $a$ we denote its degree by $\bar{a}$.

A {\it module} always means a (left) graded module.

A DG algebra $\cB=(\cB ,d_{\cB})$ is a (graded) algebra with a map
$d=d_{\cB}:\cB \to \cB$ of degree 1 such that $d^2=0$, $d(1)=0$ and
$$d(ab)=d(a)b+(-1)^{\bar{a}}ad(b).$$

Given a DG algebra $\cB$ its opposite is the DG algebra $\cB ^{op}$
which has the same differential as $\cB$ and multiplication
$$a\cdot b=(-1)^{\bar{a}\bar{b}}ba,$$
where $ba$ is the product in $\cB$. When there is a danger of
confusion of the opposite DG algebra $\cB ^{op}$ with the degree zero
part of $\cB$ we will add a comment.

We denote by $\dgalg$ the category of DG algebras.

A (left) DG module over a DG algebra $\cB$ is called a DG
$\cB$-module or, simply a $\cB$-module. A {\it right} $\cB$-module
is a DG module over $\cB ^{op}$.

If $\cB$ is a DG algebra and $M$ is a usual (not DG) module over the
algebra $\cB$, then we say that $M^{\gr}$ is a $\cB ^{\gr}$-module.

An {\it augmentation} of a DG algebra $\cB$ is a (surjective)
homomorphism of DG algebras $\cB \to k$. Its kernel is a DG ideal
(i.e. an ideal closed under the differential) of $\cB$. Denote by
$\adgalg$ the category of augmented DG algebras (morphisms commute
with the augmentation).

\begin{defi} Let $R$ be an algebra. We call $R$
{\it artinian}, if it is finite dimensional  and has a (graded)
nilpotent two-sided (maximal) ideal $m\subset R$, such that $R/m=k$.
\end{defi}

\begin{defi} Let $\cR $ be an augmented DG algebra.
We call $\cR$ {\it artinian} if $\cR$ is  artinian as an algebra and
the maximal ideal $m\subset R$ is a DG ideal, i.e. the quotient map
$R\to R/m$ is an augmentation of the DG algebra $\cR$. Note that a
homomorphism of artinian DG algebras automatically commutes with the
augmentations. Denote by $\dgart$ the category of artinian DG
algebras.
\end{defi}

\begin{defi} An artinian DG algebra $\cR$ is called positive (resp. negative) if
negative (resp. positive) degree components of $\cR$ are zero.
Denote by $\dgart _+$ and $\dgart _-$ the corresponding full
subcategories of $\dgart$. Let $\art :=\dgart _-\cap \dgart _+$ be
the full subcategory of $\dgart$ consisting of (not necessarily
commutative) artinian algebras concentrated in degree zero. Denote
by $\cart \subset \art$ the full subcategory of commutative artinian
algebras.
\end{defi}

Given a DG algebra $\cB$ one studies the category $\cB\text{-mod}$
and the corresponding homotopy and derived categories. A
homomorphism of DG algebras induces various functors between these
categories. We will recall these categories and functors in the more
general context of DG categories in the next section.

\section{DG categories}

In this section we recall some basic facts about DG categories which
will be needed in this paper. Our main references here are
\cite{BoKa}, \cite{Dr}, \cite{Ke}.

 A DG category is a $k$-linear category $\cA$ in which the sets $\Hom (A,B)$, $A,B\in Ob\cA$,
 are
 provided
with a structure of a $\bbZ$-graded $k$-module and a differential
$d:\Hom(A,B)\to \Hom (A,B)$ of degree 1, so that for every $A,B,C\in
\cA$ the composition $\Hom (A,B)\times\Hom (B,C)\to \Hom (A,C)$
comes from a morphism of complexes $\Hom (A,B)\otimes \Hom (B,C)\to
\Hom (A,C)$. The identity morphism $1_A\in \Hom (A,A)$ is closed of
degree zero.

The simplest example of a DG category is the category $DG(k)$ of
complexes of $k$-vector spaces, or DG $k$-modules.

Note also that a DG algebra is simply a DG category with one object.

Using the supercommutativity isomorphism $S\otimes T\simeq T\otimes
S$ in the category of DG $k$-modules one defines for every DG
category $\cA$ the opposite DG category $\cA ^{op}$ with $Ob\cA
^{op}=Ob\cA$, $\Hom_{\cA ^{op}}(A,B)=\Hom _{\cA}(B,A)$. We denote by $\cA
^{\gr}$ the {\it graded }
 category which is obtained from $\cA$ by forgetting the differentials on $\Hom $'s.

The tensor product of DG-categories $\cA$ and $\cB$ is defined as
follows:

(i) $Ob(\cA \otimes \cB):=Ob\cA \times Ob\cB$; for $A\in Ob\cA$ and
$B\in Ob\cB$ the corresponding object is denoted by $A\otimes B$;

(ii) $\Hom(A\otimes B,A^\prime \otimes B^\prime):=\Hom
(A,A^\prime)\otimes \Hom (B,B^\prime)$ and the composition map is
defined by $(f_1\otimes g_1)(f_2\otimes g_2):=
(-1)^{\bar{g_1}\bar{f_2}}f_1f_2\otimes g_1g_2.$

Note that the DG categories $\cA \otimes \cB$ and $\cB \otimes \cA$
are canonically isomorphic. In the above notation the isomorphism DG
functor $\phi$ is
$$\phi (A\otimes B)=(B\otimes A), \quad \phi(f\otimes g)=(-1)^{\bar{f}\bar{g}}(g\otimes f).$$

Given a DG category $\cA$ one defines the graded category
$\Ho^\bullet (\cA)$ with $Ob\Ho^\bullet (\cA)=Ob\cA$ by replacing
each $\Hom$ complex by the direct sum of its cohomology groups. We
call $\Ho^\bullet (\cA)$ the {\it graded homotopy category} of
$\cA$. Restricting ourselves to the 0-th cohomology of the $\Hom $
complexes we get the {\it homotopy category} $\Ho(\cA)$.

Two objects $A,B\in Ob\cA$ are called DG {\it isomorphic} (or,
simply, isomorphic) if there exists an invertible degree zero
morphism $f\in \Hom(A,B)$. We say that $A,B$ are {\it homotopy
equivalent} if they are isomorphic in  $\Ho(\cA)$.

A DG-functor between DG-categories $F:\cA \to \cB$ is said to be a
{\it quasi-equivalence} if $\Ho ^\bullet(F):\Ho ^\bullet(\cA)\to
\Ho ^\bullet(\cB)$ is an equivalence of graded categories. We say
that $F$ is a DG {\it equivalence} if it is fully faithful and
every object of $\cB$ is DG isomorphic to an object of $F(\cA)$.
Certainly, a DG equivalence is a quasi-equivalence. DG categories
$\cC$ and $\cD$ are called {\it quasi-equivalent} if there exist
DG categories $\cA _1,...,\cA _n$ and a chain of
quasi-equivalences
$$\cC \leftarrow \cA _1 \rightarrow ...\leftarrow \cA _n \rightarrow \cD.$$

Given DG categories $\cA$ and $\cB$ the collection of covariant DG
functors $\cA \to \cB$ is itself the collection of objects of a DG
category, which we denote by $\Fun _{\DG}(\cA ,\cB)$. Namely, let
$\Phi $ and $\Psi$ be two DG functors. Put $\Hom ^k(\Phi ,\Psi)$
equal to the set of natural transformations $t:\Phi ^{\gr} \to \Psi
^{\gr}[k]$ of graded functors from $\cA ^{\gr}$ to $\cB ^{\gr}$.
This means that for any morphism $f \in \Hom_{\cA}^s(A,B)$ one has
$$\Psi (f )\cdot t(A)=(-1)^{ks}t(B)\cdot \Phi (f).$$
On each $A\in \cA$ the differential of the transformation $t$ is
equal to $d(t(A))$ (one easily checks that this is well defined).
Thus, the closed transformations of degree 0 are the DG
transformations of DG functors. A similar definition gives us the
DG-category
 consisting of the contravariant DG functors
 $\Fun _{\DG}(\cA ^{op} ,\cB)=\Fun _{\DG}(\cA  ,\cB ^{op})$
 from $\cA$ to $\cB$.

\subsection{DG modules over DG categories}
We denote the DG category $\Fun _{\DG}(\cA ,DG(k))$ by $\cA
\text{-mod}$ and call it the
 category
of DG $\cA$-modules. There is a natural covariant DG functor $h:\cA
\to \cA ^{op}\text{-mod}$ (the Yoneda embedding) defined by
$h^A(B):=\Hom _{\cA}(B,A)$. As in the "classical" case one verifies
that the functor $h$ is fully faithful, i.e. there is a natural
isomorphism of complexes
$$\Hom _{\cA}(A,A^\prime)=\Hom_{\cA ^{op}\text{-mod}}(h^A,h^{A^\prime}).$$
Moreover, for any $M\in \cA ^{op}\text{-mod}$, $A\in \cA$
$$\Hom _{\cA ^{op}\text{-mod}}(h^A,M)=M(A).$$

The DG $\cA^{op}$-modules $h^A$, $A\in \cA$ are called {\it free}.

For $A\in \cA$ one may consider also the covariant DG functor
$h_A(B):=\Hom _{\cA}(A,B)$ and the contravariant DG functor
$h^*_A(B):=\Hom _k(h_A(B),k)$. For any $M\in \cA ^{op}\text{-mod}$ we
have
$$\Hom _{\cA ^{op}\text{-mod}}(M,h^*_A)=\Hom _k(M(A),k).$$

 A DG $\cA^{op}$-module $M$ is called acyclic, if the complex $M(A)$
is acyclic for all $A\in \cA$. Let $D(\cA ^{op})$ denote the {\it
derived category} of DG $\cA ^{op}$-modules, i.e. $D(\cA ^{op})$ is the
Verdier quotient of the homotopy category $\Ho(\cA^{op}\text{-mod})$
by the subcategory of acyclic DG-modules. This is a triangulated
category.

A DG $\cA ^{op}$-module $P$ is called h-{\it projective} if for any
acyclic DG $\cA ^{op}$-module $N$ the complex $\Hom (P,N)$ is
acyclic. A free DG module is h-projective. Denote by $\P(\cA ^{op})$
the full DG subcategory of $\cA^{op}\text{-mod}$ consisting of
h-projective DG modules.

Similarly, a  DG $\cA ^{op}$-module $I$ is called h-{\it injective}
if for any acyclic DG $\cA ^{op}$-module $N$ the complex $\Hom (N,I)$
is acyclic. For any $A\in \cA$ the DG $\cA ^{op}$-module $h^*_A$ is
h-injective.  Denote by $\cI(\cA ^{op})$ the full DG subcategory of
$\cA^{op}\text{-mod}$ consisting of h-injective DG modules.

For any DG category $\cA$ the DG categories $\cA^{op}\text{-mod}$,
$\P(\cA ^{op})$, $\cI(\cA ^{op})$ are (strongly) pre-triangulated
(\cite{Dr, BoKa}, also see subsection 3.5 below). Hence the homotopy
categories $\Ho(\cA^{op}\text{-mod})$, $\Ho(\P(\cA ^{op}))$,
$\Ho(\cI(\cA ^{op}))$ are triangulated.

The following theorem was proved in \cite{Ke}.

\begin{theo}  The inclusion functors $\P(\cA ^{op} )\hookrightarrow \cA
^{op}\text{-mod}$, $\cI(\cA ^{op})\hookrightarrow \cA ^{op}\text{-mod}$
 induce  equivalences of triangulated categories $\Ho(\P(\cA ^{op}))\simeq
 D(\cA ^{op})$ and $\Ho(\cI(\cA ^{op}))\simeq
 D(\cA ^{op})$.
 \end{theo}

Actually, it will be convenient for us to use some more precise
results from \cite{Ke}. Let us recall the relevant definitions.

\begin{defi} A DG $\cA ^{op}$-module $M$ is called relatively
projective if $M$ is a direct summand of a direct sum of DG
$\cA^{op}$-modules of the form $h^A[n]$, $A\in \cA$, $n\in \bbZ$. A DG
$\cA ^{op}$-module $P$ is said to have property (P) if it admits a
filtration
$$0=F_{-1}\subset F_0\subset F_1\subset ... P$$
such that

\noindent(F1) $\cup_iF_i=P$;

\noindent(F2) the inclusion $F_i\hookrightarrow F_{i+1}$ splits as a
morphism of graded modules;

\noindent(F3) each quotient  $F_{i+1}/F_i$ is a relatively
projective DG $\cA ^{op}$-module.
\end{defi}

\begin{defi}  A DG $\cA ^{op}$-module $M$ is called relatively
injective if $M$ is a direct summand of a direct product of DG
$\cA^{op}$-modules of the form $h_A^*[n]$, $A\in \cA$, $n\in \bbZ$. A
DG $\cA ^{op}$-module $I$ is said to have property (I) if it admits a
filtration
$$I=F_{0}\supset F_1\supset ...$$
such that

\noindent(F1') the canonical morphism
$$I\to \lim_{\leftarrow}I/F_i$$
is an isomorphism;

\noindent(F2') the inclusion $F_{i+1}\hookrightarrow F_i$ splits as
a morphism of graded modules;

\noindent(F3') each quotient  $F_{i}/F_{i+1}$ is a relatively
injective DG $\cA ^{op}$-module.
\end{defi}

\begin{theo} (\cite{Ke}) a) A DG $\cA ^{op}$-module with property (P) is
$h$-projective.

b) For any $M\in \cA ^{op}\text{-mod}$ there exists a quasi-isomorphism
$P\to M$, such that the DG $\cA ^{op}$-module $P$ has property (P).

c)  A DG $\cA ^{op}$-module with property (I) is $h$-injective.

d) For any $M\in \cA ^{op}\text{-mod}$ there exists a quasi-isomorphism
$M\to I$, such that the DG $\cA ^{op}$-module $I$ has property (I).
\end{theo}

\begin{remark} a) Assume that a DG $\cA ^{op}$-module $M$ has an
increasing filtration $M_1\subset M_2\subset ...$ such that $\cup
M_i=M$, each inclusion $M_i\hookrightarrow M_{i+1}$ splits as a
morphism of graded modules,  and each subquotient $M_{i+1}/M_i$ is
$h$-projective. Then $M$ is h-projective. b) Assume that a DG $\cA
^{op}$-module $N$ has a decreasing filtration $N=N_1\supset N_2\supset
...$ such that $\cap N_i=0$, each inclusion $N_{i+1}\hookrightarrow
N_i$ splits as a morphism of graded modules,  each subquotient
$N_i/N_{i+1}$ is h-injective (hence $N/N_i$ is h-injective for each
$i$) and the natural map
$$N\to \lim_{\leftarrow}N/N_i$$
is an isomorphism. Then $N$ is h-injective.
\end{remark}

\subsection{Some DG functors}
Let $\cB$ be a small DG category. The complex
$$\Alg _{\cB}:=\bigoplus _{A,B\in Ob \cB}\Hom(A,B)$$
has a natural structure of a DG algebra possibly without a unit. It
has the following property: every finite subset of $\Alg _{\cB}$ is
contained in $e\Alg _{\cB} e$ for some idempotent $e$ such that
$de=0$ and $\bar{e}=0$. We say that a DG module $M$ over $\Alg
_{\cB}$ is {\it quasi-unital} if every element of $M$ belongs to
$eM$ for some idempotent $e\in \Alg _{\cB}$ (which may be assumed
closed of degree $0$ without loss of generality). If $\Phi $ is a DG
$\cB$-module then
$$M_{\Phi}:=\oplus _{A\in Ob \cB}\Phi (A)$$
 is a quasi-unital
DG module over $\Alg _{\cB}$. This way we get a DG equivalence
between DG category of DG $\cB$-modules and that of quasi-unital DG
modules over $\Alg _{\cB}$.

Recall that a homomorphism of (unital) DG algebras $\phi :\cA \to
\cB$ induces functors
$$\phi _*:\cB^{op}\text{-mod}\to \cA^{op}\text{-mod},$$
$$\phi ^*:\cA^{op}\text{-mod}\to \cB^{op} \text{-mod}$$
$$\phi ^!:\cA^{op}\text{-mod}\to \cB^{op} \text{-mod}$$
where $\phi _*$ is the restriction of scalars, $\phi ^*(M)=M \otimes
_{\cA}\cB$ and $\phi ^!(M)=\Hom _{\cA ^{op}}(\cB,M)$. The DG functors
$(\phi ^*,\phi _*)$  and $(\phi _*,\phi ^!)$ are adjoint: for $M\in
\cA^{op}\text{-mod}$ and $N\in \cB^{op}\text{-mod}$ there exist
functorial isomorphisms of complexes
$$\Hom (\phi ^*M,N)=\Hom (M,\phi _*N),\quad
\Hom (\phi _*N,M)=\Hom (N,\phi ^!M). $$

This generalizes to a DG functor $F:\cA \to \cB$ between DG
categories. We obtain DG functors
$$F _*:\cB^{op}\text{-mod}\to \cA^{op}\text{-mod},$$
$$F ^*:\cA^{op}\text{-mod}\to \cB^{op}\text{-mod}.$$
$$F ^!:\cA^{op}\text{-mod}\to \cB^{op}\text{-mod}.$$

Namely, the DG functor $F$ induces a homomorphism of DG algebras
$F:\Alg _{\cA}\to \Alg _{\cB}$ and hence defines functors $F_*$,
$F^*$ between quasi-unital DG modules as above. (These functors
$F_*$ and $F^*$ are denoted in \cite{Dr} by $\Res _F$ and $\Ind _F$
respectively.) The functor $F^!$ is defined as follows: for a
quasi-unital $\Alg _{\cA}^{op}$-module $M$ put
$$F^!(M)=\Hom _{\Alg _{\cA}^{op}}(\Alg _{\cB},M)^{\qu},$$
where $N^{\qu}\subset N$ is the {\it quasi-unital} part of a $\Alg
_{\cB}^{op}$-module $N$ defined by
$$N^{\qu}:=\im (N\otimes _k  \Alg _{\cB}\to N).$$

The DG functors $(F ^*,F _*)$ and $(F_*,F^!)$ are adjoint.

\begin{lemma} Let $F:\cA \to \cB$ be a DG functor. Then

a) $F_*$ preserves acyclic DG modules;

b) $F^*$ preserves h-projective DG modules;

c) $F^!$ preserves h-injective DG modules.

\end{lemma}

\begin{proof} The first assertion is obvious and the other two
follow by adjunction.
\end{proof}

By Theorem 3.1 above the DG subcategories $\P(\cA ^{op})$ and $\cI(\cA
^{op})$ of $\cA ^{op}\text{-mod}$ allow us to define (left and right)
derived functors of DG functors $G:\cA ^{op}\text{-mod}\to \cB
^{op}\text{-mod}$ in the usual way. Namely for a DG $\cA ^{op}$-module
$M$ choose quasi-isomorphisms $P\to M$ and $M\to I$ with $P\in
\P(\cA ^{op})$ and $I\in \cI(\cA ^{op})$. Put
$$\bL G(M):=G(P),\quad \quad \bR G(M):=G(I).$$
In particular for a DG functor $F:\cA \to \cB$ we will consider
 derived functors $\bL F^*:D(\cA ^{op})\to D(\cB ^{op})$, $\bR F^!:D(\cA ^{op})\to D(\cB ^{op})$. We also
have the obvious functor $F_*:D(\cB ^{op})\to D(\cA ^{op})$. The functors
$(\bL F^*,F_*)$  and $(F_*, \bR F^!)$ are adjoint.

\begin{prop} Assume that the DG functor $F:\cA \to \cB$ is a
quasi-equivalence. Then

a) $F^*:\P(\cA ^{op})\to \P(\cB ^{op})$ is a quasi-equivalence;

b) $\bL F ^*:D(\cA ^{op})\to D(\cB ^{op})$ is an equivalence;

c) $F_*:D(\cB ^{op})\to D(\cA ^{op})$ is an equivalence.

d) $\bR F^!:D(\cA ^{op})\to D(\cB ^{op})$ is an equivalence.

e) $F^!:\cI (\cA ^{op})\to \cI (\cB ^{op})$ is a quasi-equivalence.

\end{prop}

\begin{proof} a) is proved in \cite{Ke} and it implies b) by
Theorem 3.1. c) (resp. d)) follows from b) (resp. c) by adjunction.
Finally, e) follows from d) by Theorem 3.1.
\end{proof}

Given DG $\cA ^{op}$-modules $M,N$ we denote by $\Ext ^n(M,N)$ the
group of morphisms $\Hom ^n _{D(\cA)}(M,N)$.

\subsection{DG category $\cA _{\cR}$} Let $\cR$ be a DG
algebra. We may and will consider $\cR$ as a DG category with one
object whose endomorphism DG algebra is $\cR$. We denote this DG
category again by $\cR$. Note that the DG category
$\cR^{op}\text{-mod}$ is just the category of right DG modules over the
DG algebra $\cR$.

For a DG category $\cA$ we denote the DG category $\cA \otimes \cR$
by $\cA _{\cR}$. Note that the collections of objects of $\cA$ and
$\cA _{\cR}$ are naturally identified. A homomorphism of DG algebras
$\phi :\cR\to \cQ$ induces the obvious DG functor $\phi=\id \otimes
\phi :\cA _{\cR}\to \cA _{\cQ}$ (which is the identity on objects),
whence the DG functors $\phi _*$, $ \phi ^*$, $\phi ^!$ between the
DG categories $\cA^{op} _{\cR}\text{-mod}$ and $\cA ^{op}
_{\cQ}\text{-mod}$. For  $M \in \cA_{\cR}^{op} \text{-mod}$
 we have
 $$\phi ^*(M)=M\otimes _{\cR}{\cQ}.$$
In case $\cQ^{\gr}$ is a finitely generated $\cR ^{\gr}$-module we
have
$$ \phi ^!(M)=\Hom _{\cR^{op}}(\cQ ,M).$$

In particular, if $\cR$ is  augmented then the canonical
homomorphisms of DG algebras $p:k\to \cR$ and $i:\cR \to k$ induce
functors
$$p:\cA \to \cA _{\cR},\quad i:\cA _{\cR}\to \cA,$$
such that $i\cdot p=\Id_{\cA}$.  So for $S\in \cA ^{op}\text{-mod}$ and
$T\in \cA ^{op}_{\cR}\text{-mod}$ we have

$$p^*(S)=S\otimes _k\cR, \quad i^*(T)=T\otimes _{\cR}k, \quad i^!(T)=\Hom _{\cR^{op}}(k,T).$$

For an artinian DG algebra $\cR$ we denote by $\cR^*$ the DG $\cR
^{op}$-module $\Hom _k(\cR,k)$. This is a left $\cR$-module by the
formula
$$rf(q):=(-1)^{(\bar{f}+\bar{q})\bar{r}}f(qr)$$
and a right $\cR$-module by the formula
$$fr(p):=f(rp)$$
for $r,p\in \cR$ and $f\in \cR ^*$. The augmentation map $\cR \to k$
defines the canonical (left and right) $\cR$-submodule $k\subset \cR
^*$. Moreover, the embedding $k\hookrightarrow \cR ^*$ induces an
isomorphism $k\to \Hom _{\cR}(k,\cR ^*)$.

\begin{defi} Let $\cR$ be an artinian DG algebra. A DG $\cA^{op}
_{\cR}$-module $M$ is called graded $\cR$-free (resp. graded
$\cR$-cofree) if there exists a DG $\cA ^{op}$-module $K$ such that
$M^{\gr}\simeq (K\otimes \cR)^{\gr}$ (resp. $M^{\gr}\simeq (K\otimes
\cR^*)^{\gr}$). Note that for such $M$ one may take $K=i^*M$ (resp.
$K=i^!M$).
\end{defi}

\begin{lemma} Let $\cR$ be an artinian DG algebra.

a) The full DG subcategories of DG $\cA _{\cR}^{op}$-modules consisting
of graded $\cR$-free (resp. graded $\cR$-cofree)  modules are DG
isomorphic. Namely, if $M\in \cA_{\cR}^{op}\text{-mod}$ is graded
$\cR$-free (resp. graded $\cR$-cofree) then $M\otimes _{\cR}\cR ^*$
(resp. $\Hom _{\cR^{op}}(\cR ^*,M)$) is graded $\cR$-cofree (resp.
graded $\cR$-free).

b) Let $M$ be a graded $\cR$-free module. There is a natural
isomorphism of DG $\cA ^{op}$-modules
$$i^*M\stackrel{\sim}{\to}i^!(M\otimes _{\cR}\cR ^*).$$
\end{lemma}

\begin{proof} a) If $M$ is graded $\cR$-free, then obviously $M\otimes
_{\cR}\cR ^*$ is graded $\cR$-cofree. Assume that $N$ is graded
$\cR$-cofree, i.e. $N^{\gr}=(K\otimes \cR ^*)^{\gr}$. Then
$$(\Hom
_{\cR^{op}}(\cR ^*,N))^{\gr}=(K\otimes \Hom _{\cR^{op}}(\cR
^*,\cR^*))^{\gr},$$ since $\dim _k\cR <\infty$. On the other hand
$$\Hom _{\cR^{op}}(\cR ^*,\cR^*)=\Hom _{\cR^{op}}(\cR ^*,\Hom _k(\cR ,k))=\Hom _k(\cR ^*\otimes
_{\cR}\cR,k)=\cR,$$ so $(\Hom _{\cR^{op}}(\cR ^*,N))^{\gr}=(K\otimes
\cR )^{\gr}$.

b) For an arbitrary DG $\cA _{\cR}^{op}$-module $M$ we have a natural
(closed degree zero) morphism of DG $\cA ^{op}$-modules
$$i^*M\to i^!(M\otimes _{\cR}\cR^*),\quad m\otimes 1\mapsto (1\mapsto m\otimes i),$$
where $i:\cR \to k$ is the augmentation map. If $M$ is graded
$\cR$-free this map is an isomorphism.
\end{proof}

\begin{prop} Let $\cR$ be an artinian DG algebra. Assume that a DG $\cA ^{op}_{\cR}$-module $M$
satisfies property (P) (resp. property (I)). Then $M$ is graded
$\cR$-free (resp. graded $\cR$-cofree).
\end{prop}

\begin{proof}  Notice that the collection of graded $\cR$-free
objects in $\cA ^{op}_{\cR}\text{-mod}$ is closed under taking direct
sums, direct summands (since the maximal ideal $m\subset \cR$ is
nilpotent) and direct products (since $\cR$ is finite
dimensional). Similarly for graded $\cR$-cofree objects since the
DG functors in Lemma 3.9 a) preserve direct sums and products.
Also notice that for any $A \in \cA_{\cR}$ the DG $\cA
^{op}_{\cR}$-module $h^A$ (resp. $h_A^*$) is graded $\cR$-free (resp.
graded $\cR$-cofree). Now the proposition follows since a DG $\cA
^{op}_{\cR}$-module $P$ (resp. $I$) with property (P) (resp. property
(I)) as a graded module is a direct sum of relatively projective
DG modules (resp. a direct product of relatively injective DG
modules).
\end{proof}

\begin{cor} Let $\cR$ be an artinian DG algebra. Then for any DG
$\cA^{op}_{\cR}$-module $M$ there exist quasi-isomorphisms $P\to M$ and
$M\to I$ such that $P\in \cP(\cA _{\cR}^{op})$, $I\in \cI (\cA
_{\cR}^{op})$ and $P$ is graded $\cR$-free, $I$ is graded $\cR$-cofree.
\end{cor}

\begin{proof} Indeed, this follows from Theorem 3.4 and Proposition
3.10 above.
\end{proof}

\begin{prop} Let $\cR$ be an artinian
DG algebra and $S,T\in \cA ^{op}_{\cR}\text{-mod}$ be graded $\cR$-free
(resp. graded $\cR$-cofree).

a) There is an isomorphism of graded vector spaces
$\Hom(S,T)=\Hom(i^*S,i^*T) \otimes \cR$, (resp.
$\Hom(S,T)=\Hom(i^!S,i^!T) \otimes \cR$), which is an isomorphism of
algebras if $S=T.$  In particular, the map $i ^*:\Hom (S,T)\to \Hom
(i^*S,i^*T)$ (resp. $i ^!:\Hom(S,T)\to \Hom(i^!S,i^!T)$) is
surjective.

b) The DG module $S$ has a finite filtration with subquotients
isomorphic to $i^*S$ as DG $\cA ^{op}$-modules (resp. to $i^!S$ as DG
$\cA^{op}$-modules).

c) The DG algebra $\End(S)$ has a finite filtration by DG ideals
with subquotients isomorphic to $\End (i^*S)$ (resp.
$\End(i^!S)$).

d) If $f\in \Hom (S,T)$ is a closed morphism of degree zero such
that $i^*f$ (resp. $i^!f$) is an isomorphism or a homotopy
equivalence or a quasi-isomorphism, then $f$ is also such.
\end{prop}

\begin{proof} Because of Lemma 3.9 above it suffices to
prove the proposition for graded $\cR$-free modules. So assume that
$S$, $T$ are graded $\cR$-free.

a) This holds because $\cR$ is finite dimensional.

b) We can refine the filtration of $\cR$ by powers of the maximal
ideal to get a filtration  $F_i\cR$ by ideals with 1-dimensional
subquotients (and zero differential). Then the filtration
$F_iS:=S\cdot F_i\cR$ satisfies the desired properties.

c) Again the filtration $F_i\End (S):=\End(S)\cdot F_i\cR$ has the
desired properties.

d) If $i^*f$ is an isomorphism, then $f$ is surjective by the
Nakayama lemma for $\cR$. Also $f$ is injective since $T$ is
graded $\cR$-free.

Assume that $i^*f$ is a homotopy equivalence. Let $C(f)\in \cA
_{\cR}^{op}\text{-mod}$ be the cone of $f$. (It is also graded
$\cR$-free.) Then $i^*C(f)\in \cA ^{op}\text{-mod}$ is the cone
$C(i^*f)$ of the morphism $i^*f$. By assumption the DG algebra $\End
(C(i^*f))$ is acyclic. But by part c) the complex $\End (C(f))$ has
a finite filtration with subquotients isomorphic to the complex
$\End (C(i^*f))$. Hence $\End (C(f))$ is also acyclic, i.e. the DG
module $C(f)$ is null-homotopic, i.e. $f$ is a homotopy equivalence.

Assume that $i^*f$ is a quasi-isomorphism. Then in the above
notation $C(i^*f)$ is acyclic. Since by part b) $C(f)$ has a finite
filtration with subquotients isomorphic to $C(i^*f)$, it is also
acyclic. Thus $f$ is a quasi-isomorphism.
\end{proof}

\subsection{More DG functors} So far we considered DG functors $F_*$, $F^*$,
$F^!$ between the DG categories $\cA ^{op}$-mod and $\cB ^{op}$-mod which
came from a DG functor $F:\cA \to \cB$. We will also need to
consider a different type of DG functors.

\begin{example}  For an artinian DG algebra $\cR$ and a small DG category $\cA$
we will consider two types of "restriction of scalars" DG functors
$\pi _*,\pi _!:\cA ^{op}_{\cR}\text{-mod}\to \cR ^{op}\text{-mod}$.
Namely, for $M\in \cA _{\cR}^{op}\text{-mod}$ put
$$\pi _*M:=\prod_{A\in Ob\cA _{\cR}}M(A),\quad \pi _!M:=\bigoplus_{A\in Ob\cA
_{\cR}}M(A).$$ We will also consider the two "extension of scalars"
functors $\pi ^*,\pi ^!:\cR ^{op}\text{-mod}\to \cA
^{op}_{\cR}\text{-mod}$ defined by
$$\pi ^*(N)(A):=N\otimes \bigoplus_{B\in Ob\cA }\Hom_{\cA }(A,B),
\quad \pi ^!(N)(A):=\Hom _k(\bigoplus_{B\in Ob\cA }\Hom_{\cA
}(B,A),N)$$ for  $A\in Ob\cA _{\cR}$. Notice that the DG functors
$(\pi ^*,\pi _*)$ and $(\pi _!,\pi ^!)$ are adjoint, that is for
$M\in \cA ^{op}_{\cR}\text{-mod}$ and $N\in \cR ^{op}\text{-mod}$ there is
a functorial isomorphism of complexes
$$\Hom (\pi ^*N,M)=\Hom (N,\pi_*M), \quad \Hom (\pi _!M,N)=\Hom (M,\pi^!N).$$

The DG functors $\pi^*,\pi ^!$  preserve acyclic DG modules, hence
$\pi _*$ preserves h-injectives and $\pi _!$ preserves
h-projectives.

We have the following commutative functorial diagrams
$$\begin{array}{ccc}
\cA _{\cR}^{op}\text{-mod} & \stackrel{i^*}{\longrightarrow} & \cA
^{op}\text{-mod}\\
\pi _!\downarrow & & \pi _!\downarrow \\
\cR ^{op}\text{-mod} & \stackrel{i^*}{\longrightarrow} & DG(k),
\end{array}$$
$$\begin{array}{ccc}
\cA _{\cR}^{op}\text{-mod} & \stackrel{i^!}{\longrightarrow} & \cA
^{op}\text{-mod}\\
\pi _*\downarrow & & \pi _*\downarrow \\
\cR ^{op}\text{-mod} & \stackrel{i^!}{\longrightarrow} & DG(k).
\end{array}$$
\end{example}

\begin{example} Fix $E\in \cA ^{op}\text{-mod}$ and put $\cB=\End(E)$.
Consider the DG functor
$$\Sigma =\Sigma ^E:\cB ^{op}\text{-mod}\to \cA ^{op}\text{-mod}$$
defined by $\Sigma (M)=M\otimes _{\cB}E$. Clearly, $\Sigma (\cB)=E$.
This DG functor gives rise to the functor
$$\bL \Sigma :D(\cB ^{op})\to D(\cA ^{op}),\quad \quad \bL \Sigma
(M)=M\stackrel{\bL}{\otimes }_{\cB}E.$$
 \end{example}

\subsection{Pre-triangulated DG categories} For any DG category $\cA$ there
exists a DG category $\cA ^{pre-tr}$ and a canonical full and
faithful DG functor $F:\cA \to \cA^{\pre-tr}$ (see \cite{BoKa, Dr}).
The homotopy category $\Ho (\cA^{\pre-tr})$ is canonically
triangulated. The DG category $\cA$ is called {\it pre-triangulated}
if the DG functor $F$ is a quasi-equivalence. The DG category
$\cA^{\pre-tr}$ is pre-triangulated.

Let $\cB$ be another DG category and $G:\cA \to \cB$ be a
quasi-equivalence. Then $G^{\pre-tr}:\cA ^{\pre-tr}\to \cB
^{\pre-tr}$ is also a quasi-equivalence.

The DG functor $F$ induces a DG isomorphism of DG categories
$F_*:(\cA^{\pre-tr})^{op}\text{-mod} \to \cA ^{op}\text{-mod}$. Hence the
functors $F_*:D((\cA^{\pre-tr})^{op})\to D(\cA ^{op})$ and $\bL F^*:D(\cA
^{op})\to D((\cA ^{\pre-tr})^{op})$ are equivalences. We obtain the
following corollary.

\begin{cor} Assume that a DG functor $G_1:\cA \to \cB$ induces a quasi-equivalence
$G_1^{\pre-tr}:\cA ^{\pre-tr}\to \cB ^{\pre-tr}$. Let $\cC$ be
another DG category and consider the DG functor $G:=G_1\otimes \id
:\cA \otimes \cC\to \cB \otimes \cC$.
 Then the functors $G_*,\bL G^*, \bR G^!$ between
the derived categories $D((\cA \otimes \cC)^{op})$ and $D((\cB \otimes
\cC)^{op})$ are equivalences.
\end{cor}

\begin{proof} The DG functor $G$ induces the quasi-equivalence
$G ^{\pre-tr}:(\cA \otimes \cC)^{\pre-tr}\to (\cB \otimes \cC
)^{\pre-tr}$. Hence the corollary follows from the above discussion
and Proposition 3.6.
\end{proof}

\begin{example} Suppose $\cB$ is a pre-triangulated DG category. Let
$G_1:\cA \hookrightarrow \cB$ be an embedding of a full DG
subcategory so that the triangulated category $\Ho(\cB)$ is
generated by the collection of objects $G_1(Ob\cA)$. Then the
assumptions of the previous corollary hold.
\end{example}

\subsection{A few lemmas}

\begin{lemma} Let $\cR$, $\cQ$ be DG algebras and $M$ be a DG
$\cQ\otimes \cR ^{op}$-module.

a) For any DG modules $N$, $S$ over the DG algebras $\cQ ^{op}$ and
$\cR ^{op}$ respectively there is a natural isomorphism of complexes
$$\Hom _{\cR^{op}}(N\otimes _{\cQ}M,S)\stackrel{\sim}{\to}\Hom _{\cQ ^{op}}(N,\Hom
_{\cR^{op}}(M,S)).$$

b) There is a natural quasi-isomorphism of complexes
$$\bR \Hom _{\cR^{op}}(N\stackrel{\bL}{\otimes }_{\cQ}M,S)\stackrel{\sim}{\to}\bR \Hom _{\cQ ^{op}}(N,\bR \Hom
_{\cR^{op}}(M,S)).$$

\end{lemma}

\begin{proof} a) Indeed, for $f\in \Hom _{\cR^{op}}(N\otimes _{\cQ}M,S)$
define $\alpha (f)\in \Hom _{\cQ}(N,\Hom _{\cR^{op}}(M,S))$ by the
formula $\alpha (f)(n)(m)=f(n\otimes m)$. Conversely, for $g\in \Hom
_{\cQ}(N,\Hom _{\cR^{op}}(M,S))$ define $\beta (g)\in \Hom
_{\cR^{op}}(N\otimes _{\cQ}M,S)$ by the formula $\beta(g)(n\otimes
m)=g(n)(m)$. Then $\alpha $ and $\beta$ are mutually inverse
isomorphisms of complexes.

b) Choose quasi-isomorphisms $P\to N$ and $S\to I$, where $P\in
\cP(\cQ ^{op})$ and $I\in \cI (\cR ^{op})$ and apply a).
\end{proof}

\begin{lemma} Let $\cR$ be an artinian DG algebra. Then in the DG
category $\cR ^{op}\text{-mod}$ a direct sum of copies of $\cR ^*$ is
h-injective.
\end{lemma}

\begin{proof} Let $V$ be a graded vector space, $M=V\otimes \cR
^*\in \cR ^{op}\text{-mod}$ and $C$ an acyclic DG $\cR^{op}$-module.
Notice that $M=\Hom _k(\cR ,V)$ since $\dim \cR <\infty$. Hence the
complex
$$\Hom _{\cR ^{op}}(C,M)=\Hom _{\cR ^{op}}(C,\Hom _k(\cR ,V))=\Hom
_k(C\otimes _{\cR}\cR,V)=\Hom _k(C,V)$$ is acyclic.
\end{proof}

\begin{lemma} Let $\cB$ be a DG algebra, such that $\cB ^i=0$ for
$i>0$. Then the category $D(\cB ^{op})$ has truncation functors: for
any DG $\cB$-module $M$ there exists a short exact sequence in the
abelian category $Z^0(\cB\text{-mod})$
$$\tau _{<0}M\to M\to \tau _{\geq 0}M,$$
where $H^i(\tau _{<0}M)=0$ if $i\geq 0$ and $H^i(\tau _{\geq 0}M)=0$
for $i<0$.
\end{lemma}

\begin{proof} Indeed, put $\tau _{<0}M:=\oplus _{i<0}M^i\oplus d(M^{-1})$.
\end{proof}

\begin{lemma} Let $\cB$ be a DG algebra, s.t. $\cB ^i=0$ for $i>0$
and $\dim \cB ^i<\infty$ for all $i$. Let $N$ be a DG $\cB$-module
with finite dimensional cohomology. Then there exists an
h-projective DG $\cB$-module $P$ and a quasi-isomorphism $P\to N$,
where $P$ in addition satisfies the following conditions

a) $P^i=0$ for $i>>0$,

b) $\dim P^i<\infty$ for all $i$.
\end{lemma}

\begin{proof} First assume that $N$ is concentrated in one degree,
say $N^i=0$ for $i\neq 0$. Consider $N$ as a $k$-module and put
$P_0:=\cB \otimes N$. We have a natural surjective map of DG
$\cB$-modules $\epsilon :P_0\to N$ which is also surjective on the
cohomology. Let $K:=\Ker \epsilon$. Then $K^i=0$ for $i>0$ and $\dim
K^i<\infty$ for all $i$. Consider $K$ as a DG $k$-module and put
$P_{-1}:=\cB \otimes K$. Again we have a surjective map of DG
$\cB$-modules $P_{-1}\to K$ which is surjective and surjective on
cohomology. And so on. This way we obtain an exact sequence of DG
$\cB$-modules
$$...\to P_{-1}\to P_0\stackrel{\epsilon}{\to}N\to 0,$$
where $P_{-j}^i=0$ for $i>0$ and $\dim P_{-j}^i<\infty$ for all
$j$. Let $P:=\oplus _jP_{-j}[j]$ be the "total" DG $\cB$-module of
the complex $...\to P_{-1} \to P_0\to 0$. Then $\epsilon :P\to N$
is a quasi-isomorphism. Since each DG $\cB$-module $P_{-j}$ has
the property (P), the module $P$ is h-projective by Remark 3.5a).
Also $P^i=0$ for $i>0$ and $\dim P^i<\infty $ for all $i$.

How consider the general case. Let $H^s(N)=0$ and $H^i(N)=0$ for all
$i<s$. Replacing $N$ by $\tau _{\geq s}N$ (Lemma 3.19) we may and
will assume that $N^i=0$ for $i<s$. Then $M:=(\Ker d_N)\cap N^s$ is
a DG $\cB$-submodule of $N$ which is not zero. If the embedding
$M\hookrightarrow N$ is a quasi-isomorphism, then we may replace $N$
by $M$ and so we are done by the previous argument. Otherwise we
have a short exact sequence of DG $\cB$-modules
$$o\to M\to N\to N/M\to 0$$
with $\dim H(M), \dim H(N/M)<\dim H(N)$. By the induction on $\dim
H(N)$ we may assume that the lemma holds for $M$ and $N/M$. But then
it also holds for $N$.
\end{proof}

\begin{cor} Let $\cB$ be a DG algebra, s.t. $\cB ^i=0$ for $i>0$,
$\dim \cB ^i<\infty$ for all $i$ and the algebra $H^0(\B)$ is
local. Let $N$ be a DG $\cB$-module with finite dimensional
cohomology. Then $N$ is quasi-isomorphic to a finite dimensional
DG $\cB$-module.
\end{cor}

\begin{proof} By Lemma 3.20 there exists a bounded above and locally
finite DG $\cB$-module $P$ which is quasi-isomorphic to $N$. It
remains to apply the appropriate truncation functor to $P$ (Lemma
3.19).
\end{proof}

\begin{cor} Let $\cB$ be an augmented DG algebra, s.t. $\cB ^i=0$ for
$i>0$, $\dim \cB ^i<\infty$ for all $i$ and the algebra $H^0(\cB)$
is local. Denote by $\langle k\rangle \subset D(\cB)$ the
triangulated envelope of the DG $\cB $-module $k$. Let $N$ be a DG
$\cB$-module with finite dimensional cohomology. Then $N\in
\langle k\rangle$.
\end{cor}

\begin{proof} By the previous corollary we may assume that $N$ is
finite dimensional. But then an easy applying of the Nakayama
lemma for $H^0(\cB)$ shows that $N$ has a filtration by DG
$\cB$-modules with subquotients isomorphic to $k$.
\end{proof}

\begin{lemma} Let $\cB$ and $\cC$ be DG algebras. Consider the DG
algebra $\cB \otimes \cC$ and a homomorphism of DG algebras $F:\cB
\to \cB \otimes \cC$, $F(b)=b\otimes 1$. Let $N$ be an h-projective
(resp. h-injective) DG $\cB \otimes \cC$-module. Then the DG
$\cB$-module $F_*N$ is also h-projective (resp. h-injective).
\end{lemma}

\begin{proof} The assertions follow from the fact that the DG
functor $F_*:\cB \otimes \cC\text{-mod}\to \cB\text{-mod}$ has a
left adjoint DG functor $F^*$ (resp. right adjoint DG functor $F^!$)
which preserves acyclic DG modules. Indeed,
$$F^*(M)=\cC \otimes _{k}M,\quad \quad F^!(M)=\Hom _k(\cC ,M).$$
\end{proof}

\part{Deformation functors}

\section{The homotopy deformation and co-deformation pseudo-functors}

Denote by ${\bf Gpd}$ the 2-category of groupoids.

Let $\cE$ be a category and $F,G:\cE \to {\bf Gpd}$ two
pseudo-functors. A morphism $\epsilon :F\to G$ is called full and
faithful (resp. an equivalence) if for every $X\in Ob\cE$ the
functor $\epsilon _X:F(X)\to G(X)$ is full and faithful (resp. an
equivalence). We call $F$ and $G$ equivalent if there exists an
equivalence $F\to G$.

It the rest of this paper  we will usually denote by $\cA$ a fixed
DG category and by $E$ a DG $\cA^{op}$-module.

Let us define the homotopy deformation pseudo-functor $\Def ^{\h}
(E):\dgart \to {\bf Gpd}$. This functor describes "infinitesimal"
(i.e. along artinian DG algebras) deformations of $E$ in the
homotopy category of DG $\cA^{op}$-modules.

\begin{defi} Let $\cR$ be an artinian DG algebra. An object in the
groupoid $\Def _{\cR}^{\h} (E)$ is a pair $(S,\sigma)$, where $S\in
\cA _{\cR}^{op}\text{-mod}$ and $\sigma :i^*S\to E$ is an isomorphism
of DG $\cA^{op}$-modules such that the following holds: there exists an
isomorphism of graded $\cA^{op} _{\cR}$-modules $\eta :(E\otimes
\cR)^{\gr} \to S^{\gr}$ so that the composition
$$E= i^*(E\otimes \cR)
\stackrel{i^*(\eta)}{\to} i^*S\stackrel{\sigma}{\to}E$$ is the
identity.

Given objects $(S,\sigma),(S^\prime ,\sigma ^\prime)\in \Def
_{\cR}^{\h}(E)$ a map $f:(S,\sigma)\to (S^\prime,\sigma ^\prime)$ is
an isomorphism $f:S\to S^\prime$ such that $\sigma ^\prime \cdot
i^*f=\sigma$. An allowable homotopy between maps $f,g$ is a homotopy
$h:f\to g$ such that $i^*(h)=0$. We define morphisms in $\Def
_{\cR}^{\h}(E)$ to be classes of maps modulo allowable homotopies.

Note that a homomorphism of artinian DG algebras $\phi :\cR \to \cQ$
induces the functor $\phi ^*:\Def _{\cR}^{\h}(E)\to \Def
_{\cQ}^{\h}(E)$. This defines the pseudo-functor
$$\Def {^h}(E):\dgart \to {\bf Gpd}.$$
\end{defi}

We refer to objects of $\Def _{\cR}^{\h} (E)$ as homotopy
$\cR$-deformations of $E$.

The term "homotopy" in the above definition is used to distinguish
the pseudo-functor $\Def ^{\h}$ from the pseudo-functor $\Def$ of
{\it derived deformations} (Definition 10.1). It may be justified by
the fact that $\Def ^{\h}(E)$ depends (up to equivalence) only on
the isomorphism class of $E$ in $\Ho (\cA ^{op}\text{-mod})$ (Corollary
8.4 a)).

\begin{example} We call $(p^*E,\id)\in \Def _{\cR}^{\h}(E)$ the
trivial $\cR$-deformation of $E$.
\end{example}

\begin{defi} Denote by $\Def _+^{\h}(E)$, $\Def _-^{\h}(E)$, $\Def
_0^{\h}(E)$, $\Def ^{\h}_{\cl}(E)$ the restrictions of the
pseudo-functor $\Def ^{\h}(E)$ to subcategories $\dgart _+$, $\dgart
_-$, $\art$, $\cart$ respectively.
\end{defi}

Let us give an alternative description of the same deformation
problem. We will define the homotopy {\it co-deformation}
pseudo-functor $\coDef^{\h}(E)$ and show that it is equivalent to
$\Def ^{\h}(E)$. The point is that in practice one should use $\Def
^{\h}(E)$ for a h-projective $E$ and $\coDef ^{\h}(E)$ for a
h-injective $E$ (see Section 11).

For an artinian DG algebra $\cR$ recall the $\cR ^{op}$-module $\cR ^*
=\Hom _k(\cR ,k)$.

\begin{defi} Let $\cR$ be an artinian DG algebra. An object in the groupoid
$\coDef^{\h}_{\cR}(E)$ is a pair $(T, \tau)$, where $T$ is a DG $\cA
^{op}_{\cR}$-module and $\tau :E\to i^!T$ is an isomorphism of DG
$\cA^{op}$-modules so that the following holds: there exists an
isomorphism of graded $\cA ^{op}_{\cR}$-modules $\delta :T^{\gr}\to
(E\otimes \cR ^*)^{\gr}$ such that the composition
$$E \stackrel{\tau}{\to}i^!T \stackrel{i^!(\delta)}{\to} i^!(E\otimes \cR ^*)
 =E$$ is the identity.

Given objects $(T,\tau)$ and $(T^\prime,\tau ^\prime)\in \coDef
^h_{\cR}(E)$ a map $g:(T,\tau)\to (T^\prime ,\tau ^\prime)$ is an
isomorphism $f:T\to T^\prime$ such that $i^!f \cdot \tau =\tau
^\prime$. An allowable homotopy between maps $f,g$ is a homotopy
$h:f\to g$ such that $i^!(h)=0$. We define morphisms in $\coDef
_{\cR}^{\h}(E)$ to be classes of maps modulo allowable homotopies.

Note that a homomorphism of DG algebras $\phi :\cR \to \cQ$ induces
the functor $\phi ^!:\coDef _{\cR}^{\h}(E)\to \coDef
_{\cQ}^{\h}(E)$. This defines the pseudo-functor
$$\coDef ^{\h}(E):\dgart \to {\bf Gpd}.$$
\end{defi}

We refer to objects of $\coDef _{\cR}^{\h} (E)$ as homotopy
$\cR$-co-deformations of $E$.

\begin{example} For example we can take $T=E\otimes \cR ^*$ with the
differential $d_{E,R^*}:=d_E\otimes 1+1\otimes d_{R^*}$ (and $\tau
=\id$). This we consider as the {\it trivial} $\cR$-co-deformation
of $E$.
\end{example}

\begin{defi} Denote by $\coDef _+^{\h}(E)$, $\coDef _-^{\h}(E)$, $\coDef
_0^{\h}(E)$, $\coDef _{\cl}^{\h}(E)$ the restrictions of the
pseudo-functor $\coDef ^{\h}(E)$ to subcategories $\dgart _+$,
$\dgart _-$, $\art$, $\cart$ respectively.
\end{defi}

\begin{prop} There exists a natural equivalence of pseudo-functors
$$\delta =\delta ^E:\Def ^{\h} (E)\to \coDef ^{\h}(E).$$
\end{prop}

\begin{proof} We use Lemma 3.9 above. Namely, let $S$ be an
$\cR$-deformation of $E$. Then $S\otimes _{\cR}\cR ^*$ is an
$\cR$-co-deformation of $E$. Conversely, given an
$\cR$-co-deformation $T$ of $E$ the  DG $\cA ^{op}_{\cR}$-module $\Hom
_{\cR^{op}}(\cR^*,T)$ is an $\cR$-deformation of $E$. This defines
mutually inverse equivalences $\delta _{\cR}$ and $\delta
_{\cR}^{-1}$  between the groupoids $\Def _{\cR}^{\h}(E)$ and
$\coDef ^{\h} _{\cR}(E)$, which extend to morphisms between
pseudo-functors $\Def ^{\h} (E)$ and $\coDef ^{\h}(E)$. Let us be a
little more explicit.

Let $\phi :\cR \to \cQ$ be a homomorphism of artinian DG algebras
and $S\in \Def ^{\h}(E)$. Then
$$\delta _{\cQ} \cdot \phi ^*(S)=S\otimes _{\cR}Q\otimes _{\cQ}\cQ
^*=S\otimes _{\cR}\cQ^*,\quad \quad \phi ^!\cdot \delta
_{\cR}(S)=\Hom _{\cR ^{op}}(\cQ ,S\otimes _{\cR}\cR ^*).$$ The
isomorphism $\alpha _{\phi}$ of these DG $\cA _{\cQ}^{op}$-modules is
defined by $\alpha _{\phi}(s\otimes f)(q)(r):=sf(q\phi (r))$ for
$s\in S$, $f\in \cQ^*$, $q\in \cQ$, $r\in \cR$. Given another
homomorphism $\psi :\cQ \to \cQ ^\prime$ of DG algebras one checks
the cocycle condition $\alpha _{\psi \phi}=\psi ^!(\alpha
_{\phi})\cdot \alpha _{\psi}$ (under the natural isomorphisms $(\psi
\phi)^*=\psi ^* \phi ^*$, $(\psi \phi)^!=\psi ^! \phi ^!$).
\end{proof}

\section{ Maurer-Cartan pseudo-functor}

\begin{defi} For a DG algebra $\cC$ with the differential $d$ consider
the (inhomogeneous) quadratic map
$$Q:\cC ^1 \to \cC ^2; \quad Q(\alpha )=d\alpha +\alpha ^2.$$
We denote by $MC(\cC)$ the (usual) Maurer-Cartan cone
$$MC(\cC)=\{ \alpha \in \cC ^1\vert Q(\alpha )=0\}.$$
\end{defi}

Note that $\alpha \in MC(\cC)$ is equivalent to the operator
$d+\alpha :\cC \to \cC$ having square zero. Thus the set $MC(\cC)$
describes the space of "internal" deformations of the differential
in the complex $\cC$.

\begin{defi} Let $\cB$ be a DG algebra with the differential $d$
and  a nilpotent DG ideal $\cI \subset \cB$. We define the
Maurer-Cartan groupoid $\cM \cC (\cB ,\cI)$ as follows. The set of
objects of $\cM \cC (\cB ,\cI)$ is the cone $MC(\cI)$. Maps between
objects are defined by means of the gauge group $G(\cB ,\cI):=1+\cI
^0$ ($\cI ^0$ is the degree zero component of $\cI$) acting on $\cM
\cC (\cB ,\cI)$ by the formula
$$g:\alpha \mapsto g\alpha g^{-1}+gd(g^{-1}),$$
where $g\in G(\cB  ,\cI)$, $\alpha \in MC(\cI)$. (This comes from
the conjugation action on the space of differentials $g:d+\alpha
\mapsto g(d+\alpha )g^{-1}$.) So if $g(\alpha)=\beta$, we call $g$ a
map from $\alpha $ to $\beta$. Denote by $G(\alpha ,\beta)$ the
collection of such maps. We define the set $\Hom (\alpha , \beta)$
in the category $\cM \cC (\cB ,\cI)$ to consist of homotopy classes
of maps, where the homotopy relation is defined as follows. There is
an action of the group $\cI ^{-1}$ on the set $G(\alpha ,\beta)$:
$$h:g\mapsto g+d(h)+\beta h+h\alpha,$$
for $h\in \cI ^{-1}, g\in G(\alpha ,\beta)$. We call two maps {\it
homotopic}, if they lie in the same $\cI ^{-1}$-orbit.
\end{defi}

To make the category $\cM \cC (\cB ,\cI)$ well defined we need to
prove a lemma.

\begin{lemma} Let $\alpha _1, \alpha _2, \alpha _3, \alpha
_4\in MC(\cI)$ and $g_1\in G(\alpha _1,\alpha _2)$, $g_1,g_3\in
G(\alpha _2,\alpha _3)$, $g_4\in G(\alpha _3 ,\alpha _4)$. If $g_2$
and $g_3$ are homotopic, then so are $g_2g_1$ and $g_3g_1$ (resp.
$g_4g_2$ and $g_4g_3$).
\end{lemma}

\begin{proof} Omit.
\end{proof}

Let $\cC$ be another DG algebra with a nilpotent DG ideal $\cJ
\subset \cC$. A homomorphism of DG algebras $\psi :\cB \to \cC$ such
that $\psi (\cI)\subset \cJ$ induces the functor
$$\psi ^*:\cM \cC (\cB ,\cI)\to \cM \cC (\cC ,\cJ).$$

\begin{defi} Let $\cB$ be a DG algebra and $\cR$ be an artinian DG
algebra with the maximal ideal $m\subset \cR$. Denote by $\cM\cC
_{\cR}(\cB)$ the Maurer-Cartan groupoid $\cM\cC(\cB \otimes
\cR,\cB\otimes m)$. A homomorphism of artinian DG algebras $\phi
:\cR \to \cQ$ induces the functor $\phi ^*:\cM\cC_{\cR}(\cB)\to
\cM\cC_{\cQ}(\cB)$. Thus we obtain the Maurer-Cartan pseudo-functor
$$\cM\cC(\cB):\dgart \to {\bf Gpd}.$$ We denote by
$\cM\cC _+(\cB)$, $\cM\cC _-(\cB)$, $\cM\cC _0(\cB)$, $\cM\cC
_{\cl}(\cB)$ the restrictions of the pseudo-functor $\cM \cC (\cB)$
to subcategories $\dgart _+$, $\dgart _-$, $\art $, $\cart $.
\end{defi}

\begin{remark} A homomorphism of DG algebras $\psi:\cC \to \cB$
induces a morphism of pseudo-functors
$$\psi ^*:\cM\cC (\cC)\to \cM \cC (\cB).$$
\end{remark}

\section{Description of pseudo-functors $\Def ^{\h}(E)$ and $\coDef ^{\h}(E)$}

We are going to give a description of the pseudo-functor $\Def
^{\h}$ and hence also of the pseudo-functor $\coDef ^{\h}$ via the
Maurer-Cartan pseudo-functor $\cM\cC$.

\begin{prop} Let $\cA$ be a DG category and $E\in \cA^{op}\text{-mod}$.
Denote by $\cB$ the DG algebra $\End(E)$. Then there exists an
equivalence of pseudo-functors $\theta =\theta ^E: \cM\cC(\cB)\to
\Def^{\h}(E)$. (Hence also $\cM\cC(\cB)$ and $\coDef^{\h}(E)$ are
equivalent.)
\end{prop}

\begin{proof} Fix an artinian DG algebra $\cR$ with the maximal ideal $m$.
Let us define an equivalence of groupoids
$$\theta _{\cR}:\cM\cC_{\cR}(\cB)\to \Def ^{\h}_{\cR}(E).$$

Denote by $S_0=p^*E\in \cA _{\cR}^{op}\text{-mod}$ the trivial
$\cR$-deformation of $E$ with the differential $d_{E,\cR}=d_E\otimes
1+1\otimes d_{\cR}$. There is a natural isomorphism of DG algebras
$\End(S_0)=\cB \otimes \cR$.

Let $\alpha \in \cM\cC(\cB\otimes m)=\cM\cC_{\cR}(\cB)$. Then in
particular $\alpha \in \End ^1(S_0)$. Hence
$d_{\alpha}:=d_{E,\cR}+\alpha$ is an endomorphism of degree 1 of the
graded module $S_0^{\gr}$. The Maurer-Cartan condition on $\alpha$
is equivalent to $d_{\alpha}^2=0$. Thus we obtain an object
$S_{\alpha}\in \cA^{op}_{\cR}\text{-mod}$. Clearly $i^*S_{\alpha}=E$,
so that
$$\theta _{\cR}(\alpha):=(S_{\alpha},\id)\in \Def _{\cR}^{\h}(E).$$

One checks directly that this map on objects extends naturally to a
functor $\theta _{\cR}:\cM\cC_{\cR}(\cB)\to \Def ^{\h}_{\cR}(E)$.
Indeed, maps  between Maurer-Cartan objects induce isomorphisms of
the corresponding deformations; also homotopies between such maps
become allowable homotopies between the corresponding isomorphisms.

It is clear that the functors $\theta _{\cR}$ are compatible with
the functors $\phi ^*$ induced by morphisms of DG algebras $\phi
:\cR \to \cQ$. So we obtain a morphism of pseudo-functors
$$\theta :\cM\cC(\cB)\to \Def^{\h}(E).$$

It suffices to prove that $\theta _{\cR}$ is an equivalence for each
$\cR$.

\medskip

\noindent{\bf Surjective.} Let $(T,\tau )\in \Def ^{h}_{\cR}(E)$. We
may and will assume that $T^{gr}=S_0^{gr}$ and $\tau =\id$. Then
$\alpha _T:=d_T-d_{\cR,E}\in \End ^1(S_0)=(\cB\otimes \cR)^1$  is an
element in $MC(\cB \otimes \cR)$. Since $i^*\alpha _T=0$ it follows
that $\alpha _T\in \cM\cC _{\cR}(\cB)$. Thus $(T,\tau )=\theta
_{\cR}(\alpha _T)$.

\medskip

\noindent{\bf Full.} Let $\alpha, \beta \in \cM \cC _{\cR}(\cB)$.
An isomorphism between the corresponding objects $\theta
_{\cR}(\alpha)$ and $\theta _{\cR}(\beta)$ is defined by an
element $f\in \End (S_0)=(\cB \otimes \cR)$ of degree zero. The
condition $i^*f=\id _Z$ means that $f\in 1+(\cB\otimes m)^0$. Thus
$f\in G(\alpha ,\beta)$.

\medskip

\noindent{\bf Faithful.} Let $\alpha, \beta \in \cM \cC _{\cR}(\cB)$
and $f,g\in G(\alpha ,\beta)$. One checks directly that $f$ and $g$
are homotopic (i.e. define the same morphism in $\cM\cC_{\cR}(\cB)$)
if and only if there exists an allowable homotopy between $\theta
_{\cR}(f)$ and $\theta _{\cR}(g)$. This proves the proposition.
\end{proof}

\begin{cor} For $E\in \cA^{op}\text{-mod}$ the pseudo-functors $\Def
^{\h}(E)$ and $\coDef ^{\h}(E)$ depend (up to equivalence) only on
the DG algebra $\End(E)$.
\end{cor}

We will prove a stronger result in Corollary 8.2 below.

\begin{example} Let $E\in \cA ^{op}\text{-mod}$ and denote $\cB
=\End(E)$. Consider $\cB$ as a (free) right $\cB$-module, i.e.
$\cB\in \cB ^{op}\text{-mod}$. Then $\Def ^h(\cB)\simeq \Def ^h(E)$
($\simeq \coDef ^h(\cB)\simeq \coDef ^h(E)$) because $\End
(\cB)=\End (E)=\cB$. We will describe this equivalence directly in
Section 9 below.
\end{example}

\section{Obstruction Theory}

It is convenient to describe the obstruction theory for our
(equivalent) deformation pseudo-functors $\Def ^h$ and $\coDef ^h$
using the Maurer-Cartan pseudo-functor $\cM\cC(\cB)$ for a fixed  DG
algebra $\cB$.

Let $\cR$ be an artinian DG algebra with a maximal ideal $m$, such
that $m^{n+1}=0$. Put $I=m^n$, $\overline{\cR}=\cR/I$ and $\pi :\cR
\to \overline{\cR}$ the projection morphism. We have $mI=Im=0$.

Note that the kernel of the homomorphism $1 \otimes \pi:\cB \otimes
\cR\to \cB \otimes \overline{\cR }$ is the (DG) ideal $\cB\otimes
I$.  The next proposition describes the obstruction theory for
lifting objects and morphisms along the functor
$$\pi^*:\cM\cC _{\cR}(\cB)\to \cM \cC _{\overline{\cR}}(\cB).$$
It is close to \cite{GoMi}. Note however a
difference in part 3) and part 4) since we do not assume that out DG algebras
live in nonnegative dimensions (and of course we work with DG
algebras and not with DG Lie algebras).

\begin{prop}
 1). There exists a map $o_2:Ob\cM\cC _{\overline{\cR}}(\cB)\to
 H^2(\cB\otimes I)$ such that $\alpha \in Ob\cM\cC
 _{\overline{\cR}}(\cB)$ is in the image of $\pi ^*$ if and only
 if $o_2(\alpha)=0$. Furthermore if $\alpha ,\beta \in Ob\cM\cC
 _{\overline{\cR}}(\cB)$ are isomorphic, then $o_2(\alpha)=0$ if
 and only if $o_2(\beta)=0$.

 2). Let $\xi \in Ob\cM\cC_{\overline{\cR}}(\cB)$. Assume that the
 fiber $(\pi ^*)^{-1}(\xi)$ is not empty. Then there exists a simply
 transitive action of the group $Z^1(\cB\otimes I)$ on the set
 $Ob(\pi ^*)^{-1}(\xi)$. Moreover the composition of the difference
 map
 $$Ob(\pi ^*) ^{-1}(\xi)\times Ob(\pi ^*) ^{-1}(\xi)\to Z^1(\cB\otimes I)$$
 with the projection
 $$Z^1(\cB\otimes I)\to H^1(\cB\otimes I)$$
 which we denote by
 $$o_1:Ob(\pi ^*)^{-1}(\xi)\times Ob(\pi ^*)^{-1}(\xi)\to H^1(\cB\otimes
 I)$$ has the following property: for $\alpha ,\beta \in Ob(\pi
 ^*)^{-1}(\xi)$ there exists a morphism $\gamma :\alpha \to \beta$
 s.t. $\pi ^*(\gamma)=\id _{\xi}$
 if and only if $o_1(\alpha ,\beta)=0$.

 3). Let $\tilde{\alpha },\tilde{\beta}\in Ob\cM\cC _{\cR}(\cB)$ be
 isomorphic objects and let $f:\alpha \to \beta$ be a morphism
 from $\alpha =\pi ^*(\tilde{\alpha})$ to $\beta =\pi
 ^*(\tilde{\beta})$. Then there is a transitive action of
 the group $H^0(\cB\otimes I)$ on the set $(\pi ^*)^{-1}(f)$ of
 morphisms $\tilde{f}:\tilde{\alpha}\to \tilde{\beta}$ such that
 $\pi ^*(\tilde{f})=f$.

 4). In the notation of 3) suppose that the fiber $(\pi^*)^{-1}(f)$
 is non-empty. Then the kernel of the above action coincides with
 the kernel of the map \begin{equation}\label{kernel} H^0(\cB\otimes I)\to H^0(\cB\otimes m, d^{\alpha,\beta}),\end{equation} where $d^{\alpha,\beta}$ is a differential on
 the graded vector space $\cB\otimes m$ given by the formula
 $$d^{\alpha,\beta}(x)=dx+\beta x-(-1)^{\bar{x}}x\alpha.$$
 In particular the difference map
 $$o_0:(\pi ^*)^{-1}(f)\times (\pi ^*)^{-1}(f)\to \im(H^0(\cB\otimes I)\to H^0(\cB\otimes m,d^{\alpha,\beta}))$$
 has the property: if $\tilde{f},\tilde{f}^\prime\in (\pi
  ^*)^{-1}(f)$, then $\tilde{f}=\tilde{f}^\prime$ if and only if
 $o_0(\tilde{f},\tilde{f}^\prime)=0$.
 \end{prop}

 \begin{proof} 1) Let $\alpha \in Ob\cM\cC
 _{\overline{\cR}}(\cB)=MC(\cB\otimes (m/I))$. Choose
 $\tilde{\alpha}\in (\cB\otimes m)^1$ such that $\pi
 (\tilde{\alpha})=\alpha$. Consider the element
 $$Q(\tilde{\alpha})=d\tilde{\alpha}+\tilde{\alpha}^2\in
 (\cB\otimes m)^2.$$
 Since $Q(\alpha)=0$ we have
 $Q(\tilde{\alpha})\in (\cB\otimes I)^2$. We claim that
 $dQ(\tilde{\alpha})=0$. Indeed,
 $$dQ(\tilde{\alpha})=d(\tilde{\alpha}^2)=d(\tilde{\alpha})\tilde{\alpha}-\tilde{\alpha}d(\tilde{\alpha}).$$
 We have $d(\tilde{\alpha})\equiv \tilde{\alpha}^2(mod(\cB\otimes
 I)).$ Hence
 $dQ(\tilde{\alpha})=-\tilde{\alpha}^3+\tilde{\alpha}^3=0$ (since
 $I\cdot m=0$).

 Furthermore suppose that $\tilde{\alpha}^\prime\in
 (\cB \otimes m)^1$ is another lift of $\alpha$, i.e.
 $\tilde{\alpha}^\prime -\tilde{\alpha}\in (\cB \otimes I)^1$. Then
 $$Q(\tilde{\alpha}^\prime)-Q(\tilde{\alpha})=d(\tilde{\alpha}^\prime-\tilde{\alpha})+
 (\tilde{\alpha}^\prime -\tilde{\alpha})(\tilde{\alpha}^\prime
 +\tilde{\alpha})=d(\tilde{\alpha}^\prime -\tilde{\alpha}).$$
Thus the cohomology class of the cocycle $Q(\tilde{\alpha})$ is
independent of the lift $\tilde{\alpha}$. We denote this class by
$o_2(\alpha)\in H^2(\cB\otimes I)$.

If $\alpha =\pi ^*(\tilde{\alpha})$ for some $\tilde{\alpha}\in
Ob\cM\cC _{\cR}(\cB)$, then clearly $o_2(\alpha)=0$. Conversely,
suppose $o_2(\alpha)=0$ and let $\tilde{\alpha}$ be as above. Then
$dQ(\tilde{\alpha})=d\tau$ for some $\tau \in (\cB\otimes I)^1$. Put
$\tilde{\alpha}^\prime=\tilde{\alpha}-\tau$. Then
$$Q(\tilde{\alpha}^\prime)=d\tilde{\alpha}-d\tau
+\tilde{\alpha}^2-\tilde{\alpha}\tau -\tau \tilde{\alpha}+\tau
^2=Q(\tilde{\alpha})-d\tau=0.$$

Let us prove the last assertion in 1). Assume that $\pi
^*(\tilde{\alpha})=\alpha$ and $\beta =g(\alpha)$ for some $g\in
1+(\cB\otimes m/I)^0$. Choose a lift $\tilde{g}\in 1+(\cB\otimes
m)^0$ of $g$ and put $\tilde{\beta}:=\tilde{g}(\tilde{\alpha})$.
Then $\pi ^*(\tilde{\beta})=\beta$. This proves 1).

2). Let $\alpha \in Ob(\pi ^*)^{-1}(\xi)$ and $\eta \in
Z^1(\cB\otimes I)$. Then
$$Q(\alpha +\eta)=d\alpha +d\eta +\alpha
^2 +\alpha \eta +\eta \alpha +\eta ^2=Q(\alpha)+d\eta =0.$$ So
$\alpha +\eta \in Ob(\pi ^*)^{-1}(\xi)$. This defines the action of
the group $Z^1(\cB \otimes I)$ on the set $Ob(\pi ^*)^{-1}(\xi)$.

Let $\alpha ,\beta \in Ob(\pi ^*)^{-1}(\xi)$. Then $\alpha -\beta
\in (\cB \otimes I)^1$ and
$$d(\alpha -\beta)=d\alpha -d\beta +\beta (\alpha -\beta)+(\alpha
-\beta)\beta +(\alpha -\beta )^2=Q(\alpha )-Q(\beta )=0.$$ Thus
$Z^1(\cB\otimes I)$ acts simply transitively on $Ob(\pi
^*)^{-1}(\xi)$. Now let $o_1(\alpha ,\beta)\in H^1(\cB\otimes I)$ be
the cohomology class of $\alpha -\beta$. We claim that there exists
a morphism $\gamma :\alpha \to \beta$ covering $\id_{\xi}$ if and
only if $o_1(\alpha ,\beta)=0$.

Indeed, let $\gamma$ be such a morphism. Then by definition the
morphisms $\pi ^*(\gamma)$ and $\id_{\xi}$ are homotopic. That is
there exists $h\in (\cB \otimes (m/I))^{-1}$ such that
$$\id_{\xi}=\pi ^*(\gamma)+d(h)+\xi h+h\xi.$$
Choose a lifting $\tilde{h}\in (\cB \otimes m)^{-1}$ on $h$ and
replace the morphism $\gamma$ by the homotopical one
$$\delta =\gamma +d(\tilde{h})+\beta \tilde{h}+\tilde{h}\alpha.$$
Thus $\delta =1+u$, where $u\in (\cB \otimes I)^0$. But then
$$\beta =\delta \alpha \delta ^{-1}+\delta d(\delta ^{-1})=\alpha
-du,$$ so that $o_1(\alpha ,\beta)=0$.

Conversely, let $\alpha -\beta=du$ for some $u\in (\cB\otimes I)^0$.
Then $\delta =1+u$ is a morphism from $\alpha $ to $\beta$ and $\pi
^*(\delta)=\id _{\xi}$. This proves 2).

3). Let us define the action of the group $Z^0(\cB\otimes I)$ on the
set $(\pi ^*)^{-1}(f)$. Let $\tilde{f}:\tilde{\alpha}\to
\tilde{\beta}$ be a lift of $f$, and $v\in Z^0(\cB\otimes I)$. Then
$\tilde{f}+v$ also belongs to $(\pi ^*)^{-1}(f)$. If $v=du$ for
$u\in (\cB\otimes I )^{-1}$, then
$$\tilde{f}+v=\tilde{f}+du+\tilde{\beta} u+u \tilde{\alpha}$$
and hence morphisms $\tilde{f}$ and $\tilde{f}+v$ are homotopic.
This induces the action of $H^0(\cB\otimes I)$ on the set $(\pi
^*)^{-1}(f)$.

To show that this action is transitive let $\tilde{f}^\prime
:\tilde{\alpha }\to \tilde{\beta}$ be another morphism in $(\pi
^*)^{-1}(f)$. This means by definition that there exists $h\in
(\cB\otimes (m/I))^{-1}$ such that
$$f=\pi ^*(\tilde{f}^\prime)+dh+\beta h+h\alpha.$$
Choose a lifting $\tilde{h}\in (\cB \otimes m)^{-1}$ of $h$ and
replace $\tilde{f}^\prime$ by the homotopical morphism
$$\tilde{g}=\tilde{f}^\prime+d\tilde{h}+\tilde{\beta}\tilde{h}+
\tilde{h}\tilde{\alpha}.$$ Then $\tilde{g}=\tilde{f}+v$ for $v\in
(\cB \otimes I)^0$. Since $\tilde{f}, \tilde{g}:\tilde{\alpha}\to
\tilde{\beta}$ we must have that $v\in Z^0(\cB \otimes I)$. This
shows the transitivity and proves 3).

4). Suppose that for some $v\in Z^0(\cB\otimes I)$ and for some
$\tilde{f}\in (\pi^*)^{-1}(f)$ we have that
$\tilde{f}+v=\tilde{f}$. This means, by definition, that there
exists an element $h\in (\cB\otimes m)^{-1}$ such that
$d^{\alpha,\beta}(h)=v$. In other words, the class $[v]\in
H^0(\cB\otimes I)$ lies in the kernel of the map (\ref{kernel}).
This proves 4).
\end{proof}

\section{Invariance theorem and its implications}

\begin{theo} Let $\phi :\cB
\to \cC$ be a quasi-isomorphism of  DG algebras. Then the induced
morphism of pseudo-functors $$\phi ^*:\cM\cC(\cB )\to \cM\cC(\cC)$$
is an equivalence.
\end{theo}

\begin{proof}
The proof is almost the same as that of Theorem 2.4 in \cite{GoMi}.
We present it for reader's convenience and also because of the
slight difference in language: in \cite{GoMi} they work with DG Lie
algebras as opposed to DG algebras.

Fix an artinian DG algebra $\cR$ with the maximal ideal $m\subset
\cR$, such that $m^{n+1}=0$. We prove that
$$\phi ^*:\cM\cC _{\cR}(\cB)\to \cM\cC _{\cR}(\cC)$$
is an equivalence by induction on $n$. If $n=o$, then both groupoids
contain one object and one morphism, so are equivalent. Let $n>0$.
Put $I=m^n$ with the projection $\pi :\cR \to \cR
/I=\overline{\cR}$. We have the commutative functorial diagram
$$\begin{array}{ccc}
\cM\cC _{\cR}(\cB) & \stackrel{\phi ^*}{\rightarrow} & \cM\cC
_{\cR}(\cC)\\
\pi ^*\downarrow & & \downarrow \pi ^*\\
\cM \cC _{\overline{\cR}}(\cB) & \stackrel{\phi ^*}{\rightarrow} &
\cM \cC _{\overline{\cR}}(\cC).
\end{array}$$
By induction we may assume that the bottom functor is an
equivalence. To prove the same about the top one we need to analyze
the fibers of the functor $\pi ^*$. This has been done by the
obstruction theory.

We will prove that the functor
$$\phi ^*:\cM\cC _{\cR}(\cB)\to \cM\cC _{\cR}(\cC)$$
is surjective on the isomorphism classes of objects, is full and is
faithful.

\medskip

\noindent{\bf Surjective on isomorphism classes.} Let $\beta \in
Ob\cM\cC _{\cR}(\cC)$. Then $\pi ^*\beta \in Ob\cM\cC
_{\overline{\cR} }(\cC)$. By the induction hypothesis there exists
$\alpha ^\prime \in Ob\cM\cC _{\overline{\cR}}(\cC)$ and an
isomorphism $g: \phi ^*\alpha ^\prime \to \pi ^* \beta$. Now
$$H^2(\phi)o_2(\alpha ^\prime)=o_2(\phi ^*\alpha ^\prime )=
o_2(\pi ^*\beta )=0.$$ Hence $o_2(\alpha ^\prime)=0$, so there
exists $\tilde{\alpha }\in Ob\cM\cC _{\cR}(\cB)$ such that $\pi
^*\tilde{\alpha}=\alpha ^\prime$, and hence
$$\phi ^*\pi ^*\tilde{\alpha}=\pi ^*\phi ^*\tilde{\alpha}=\phi
^*\alpha ^\prime.$$

Choose a lift $\tilde{g}\in 1+(\cC\otimes m)^0$ of $g$ and put
$\tilde{\beta}=\tilde{g}^{-1}(\beta)$. Then
$$\pi ^*(\tilde{\beta})=\pi ^*(\tilde{g}^{-1}(\beta))=g^{-1}\pi
^*\beta=\phi ^*\alpha ^\prime.$$

The obstruction to the existence of an isomorphism $\phi
^*\tilde{\alpha}\to \tilde{\beta}$ covering $\id_{\pi^*(\alpha
^\prime)}$ is an element $o_1(\phi ^*(\tilde{\alpha}),
\tilde{\beta})\in H^1(\cC\otimes I)$. Since $H^1(\phi)$ is
surjective there exists a cocycle $u\in Z^1(\cB\otimes I)$ such that
$H^1(\phi)[u]=o_1(\phi ^*(\tilde{\alpha}),\tilde{\beta})$. Put
$\alpha =\tilde{\alpha}-u\in Ob\cM\cC _{\cR}(\cB)$. Then
$$\begin{array}{rcl}o_1(\phi^*\alpha ,\tilde{\beta}) & = &
o_1(\phi ^*\alpha ,\phi ^*\tilde{\alpha})+o_1(\phi
^*\tilde{\alpha},\tilde{\beta})\\
& = & H^1(\phi)o_1(\alpha
,\tilde{\alpha})+o_1(\phi^*\tilde{\alpha},\tilde{\beta})\\
& = & -H^1(\phi)[u]+o_1(\phi ^*\tilde{\alpha},\beta)=0
\end{array}
$$
This proves the surjectivity of $\phi ^*$ on isomorphism classes.

\medskip

\noindent{\bf Full.} Let $f:\phi ^*\alpha _1\to \phi ^*\alpha _2$ be
a morphism in $\cM\cC _{\cR}(\cC)$. Then $\pi ^*f$ is a morphism in
$\cM\cC _{\overline{\cR}}(\cC)$:
$$\pi ^*(f):\phi^*\pi^*\alpha _1\to \phi ^*\pi ^*\alpha _2.$$

By induction hypothesis there exists $g:\pi ^*\alpha _1\to \pi
^*\alpha _2$ such that $\phi ^*(g)=\pi^*(f)$. Let $\tilde{g}\in
1+(\cC \otimes m)^0$ be any lift of $g$. Then $\pi
^*(\tilde{g}\alpha _1)=\pi ^*\alpha _2$. The obstruction to the
existence of a morphism $\gamma :\tilde{g}\alpha _1\to \alpha _2$
covering $\id_{\pi ^*\alpha _2}$ is an element $o_1(\tilde{g}\alpha
_1,\alpha _2)\in H^1(\cB\otimes I)$. By assumption $H^1(\phi)$ is an
isomorphism and we know that
$$H^1(\phi)(o_1(\tilde{g}\alpha
_1,\alpha _2))=o_1(\phi^*\tilde{g}\alpha _1,\phi ^*\alpha _2)=0,$$
since the morphism $f\cdot (\phi ^*\tilde{g})^{-1}$ is covering
the identity morphism $\id_{\pi ^*\phi ^*\alpha _2}$. Thus
$o_1(\tilde{g}\alpha _1,\alpha _2)=0$ and $\gamma $ exists. Then
$\gamma \cdot \tilde{g}:\alpha _1\to \alpha _2$ is covering $g:\pi
^*\alpha _1\to \pi^*\alpha _2$. Hence both morphisms $\phi
^*(\gamma \cdot \tilde{g})$ and $f$ are covering $\pi ^*(f)$. The
obstruction to their equality is an element $o_0(\phi ^*(\gamma
\cdot \tilde{g}),f)\in \im(H^0(\cC \otimes I)\to H^0(\cC\otimes
m))$. Let $v\in H^0(\cC\otimes I)$ be a representative of this
element and $u\in Z^0(\cB\otimes I)$ be a representative of the
inverse image of $v$ under $H^0(\phi)$. Then $\phi ^*(\gamma \cdot
\tilde{g}+u)=f$.

\medskip

\noindent{\bf Faithful.} Let $\gamma _1,\gamma _2:\alpha _1\to
\alpha _2$ be morphisms in $\cM\cC _{\cR}(\cB)$ with $\phi^*\gamma
_1=\phi^*\gamma _2$. Then $\phi ^*\pi ^*\gamma _1=\phi ^*\pi
^*\gamma _2$. By the induction hypothesis $\pi ^*\gamma _1=\pi
^*\gamma _2$, so the obstruction $o_0(\gamma _1,\gamma _2)\in
\im(H^0(\cB\otimes I)\to H^0(\cB\otimes m,d^{\alpha_1,\alpha_2}))$
is defined. Now the image of $o_0(\gamma _1,\gamma _2)$ under the
map\begin{equation}\label{iso}\im(H^0(\cB\otimes I)\to
H^0(\cB\otimes m,d^{\alpha_1,\alpha_2})) \to \im(H^0(\cC\otimes
I)\to H^0(\cC\otimes
m,d^{\phi^*\alpha_1,\phi^*\alpha_2}))\end{equation} equals to
$o_0(\phi^*\gamma _1,\phi^*\gamma _2)=0$. So it remains to prove
that the map (\ref{iso}) is an isomorphism. Clearly, it is
sufficient to prove that the morphism of complexes
$$\phi_{\cR}^{\alpha_1,\alpha_2}:(\cB\otimes m,d^{\alpha_1,\alpha_2}))\to (\cC\otimes
m,d^{\phi^*\alpha_1,\phi^*\alpha_2}))$$ is a quasi-isomorphism.
Note that these complexes have finite filtrations by subcomplexes
$\cB\otimes m^i$ and $\cC\otimes m^i$ respectively. The morphism
$\phi_{\cR}^{\alpha_1,\alpha_2}$ is compatible with these
filtrations and induces quasi-isomorphisms on the subquotients.
Hence $\phi_{\cR}^{\alpha_1,\alpha_2}$ is a quasi-isomorphism.
This proves the theorem.
\end{proof}

\begin{cor} The homotopy (co-) deformation pseudo-functor of $E\in \cA
^{op}\text{-mod}$ depends (up to equivalence) only on the
quasi-isomorphism class of the DG algebra $\End (E)$.
\end{cor}

\begin{proof} This follows from Theorem 8.1 and Proposition 6.1.
\end{proof}

The next proposition provides two examples of this situation. It was
communicated to us by Bernhard Keller.

\begin{prop} (Keller)  a) Assume that $E^\prime \in \cA
^{op}\text{-mod}$ is homotopy equivalent to $E$. Then the DG algebras
$\End (E)$ and $\End(E^\prime )$ are canonically quasi-isomorphic.

b) Let $P\in \cP(\cA ^{op})$ and $I\in \cI(\cA ^{op})$  be
quasi-isomorphic. Then the DG algebras $\End(P)$ and $\End(I)$ are
canonically quasi-isomorphic.
\end{prop}

\begin{proof} a) Let $g:E\to E^\prime$ be a homotopy equivalence.
Consider its cone $C(g)\in \cA ^{op}\text{-mod}$. Let $\cC \subset
\End (C(g))$ be the DG subalgebra consisting of endomorphisms
which leave $E^\prime $ stable. There are natural projections
$p:\cC\to \End(E^\prime)$ and $q:\cC\to \End (E)$. We claim that
$p$ and $q$ are quasi-isomorphisms. Indeed, $\Ker(p)$ (resp. $\Ker
(q)$) is the complex $\Hom (E[1],C(g))$ (resp. $\Hom
(C(g),E^\prime)$). These complexes are acyclic, since $g$ is a
homotopy equivalence.

b) The proof is similar. Let $f:P\to I$ be a quasi-isomorphism.
Then the cone $C(f)$ is acyclic. We consider the DG subalgebra
$\D\subset \End (C(f))$ which leaves $I$ stable. Then $\D$ is
quasi-isomorphic to $\End(I)$ and $\End(P)$ because the complexes
$\Hom (P[1],C(f))$ and $\Hom (C(f),I)$ are acyclic.
\end{proof}

\begin{cor} a) If DG $\cA ^{op}$-modules $E$ and $E^\prime$ are
homotopy equivalent then the pseudo-functors $\Def ^{\h}(E)$,
$\coDef ^{\h}(E)$, $\Def ^{\h}(E^\prime)$, $\coDef ^{\h}(E^\prime)$
are canonically equivalent.

b) Let $P\to  I$ be a quasi-isomorphism between $P\in \cP(\cA ^{op})$
and $I\in \cI(\cA ^{op})$. Then the pseudo-functors $\Def ^{\h}(P)$,
$\coDef ^{\h}(P)$, $\Def ^{\h}(I)$, $\coDef ^{\h}(I)$ are
canonically equivalent.
\end{cor}

\begin{proof} Indeed, this follows from Proposition 8.3 and
Corollary 8.2.
\end{proof}

Actually, one can prove a more precise statement.

\begin{prop} Fix an artinian DG algebra $\cR$.

a) Let $g:E\to E^\prime $ be a homotopy equivalence of DG $\cA
^{op}$-modules. Assume that $(V,\id)\in \Def ^{\h}_{\cR}(E)$ and
$(V^\prime,\id)\in \Def ^{\h}_{\cR}(E^\prime)$ are objects that
correspond to each other via the equivalence $\Def
^{\h}_{\cR}(E)\simeq \Def ^{\h}_{\cR}(E^\prime)$ of Corollary 8.4.
Then there exists a homotopy equivalence $\tilde{g}:V\to V^\prime$
which extends $g$, i.e. $i^*\tilde{g}=g$. Similarly for the objects
of $\coDef ^{\h}_{\cR}$ with $i^!$ instead of $i^*$.

b) Let $f:P\to I$ be a quasi-isomorphism with $P\in \cP(\cA ^{op})$,
$I\in \cI(\cA ^{op})$.  Assume that $(S,\id)\in \Def ^{\h}_{\cR}(P)$
and $(T,\id)\in \Def ^{\h}_{\cR}(I)$ are objects that correspond to
each other via the equivalence $\Def ^{\h}_{\cR}(P)\simeq \Def
^{\h}_{\cR}(I)$ of Corollary 8.4. Then there exists a
quasi-isomorphism $\tilde{f}:S\to T$ which extends $f$, i.e.
$i^*\tilde{f}=f$. Similarly for the objects of $\coDef^{\h}_{\cR}$
with $i^!$ instead of $i^*$.
\end{prop}

\begin{proof} a)
 Consider the DG algebra
 $$\cC   \subset \End(C(g))$$ as in the proof of Proposition 8.3.
 We proved there that the natural projections $\End
 (E)\leftarrow \cC \rightarrow \End(E^\prime)$ are
 quasi-isomorphisms. Hence the induced functors between groupoids
 $\cM \cC_{\cR}(\End (E))\leftarrow \cM \cC _{\cR}(\cC)\rightarrow \cM \cC_{\cR}(\End
 (E^\prime))$ are equivalences by Theorem 8.1. Using Proposition 6.1
 we may and will assume that deformations $(V, \id)$, $(V^\prime
 ,\id )$ correspond to elements $\alpha _E\in \cM \cC_{\cR}(\End
 (E))$, $\alpha _{E^\prime}\in \cM \cC_{\cR}(\End (E^\prime))$ which
 come from the same element $\alpha \in \cM \cC_{\cR}(\cC)$.

 Consider the DG modules $E\otimes \cR$, $E^\prime \otimes \cR$
 with the differentials $d_E\otimes 1+1\otimes d_{\cR}$ and
 $d_{E^\prime}\otimes 1+1\otimes d_{\cR}$ respectively and the morphism $g\otimes 1:E\otimes \cR \to
 E^\prime \otimes \cR$. Then
$$\cC\otimes \cR =\left(  \begin{array}{cc}
\End(E^\prime\otimes \cR) & \Hom (E[1]\otimes \cR,E^\prime \otimes \cR) \\
0 & \End(E\otimes \cR)
\end{array}  \right)
  \subset \End(C(g\otimes 1)),$$
  and
 $$\alpha =\left( \begin{array}{cc}
 \alpha _{E^\prime} & t\\
 0                  & \alpha _{E}
 \end{array}\right).$$

Recall that the differential in the DG module $C(g\otimes 1)$
 is of the form $(d_{E^\prime}\otimes 1, d_E[1]\otimes 1+g[1]\otimes
 1)$. The element $\alpha$ defines a new differential $d_{\alpha}$ on
$C(g\otimes 1)$ which is $(d_{E^\prime}\otimes 1+\alpha _{E^\prime},
(d_E[1]\otimes 1+\alpha _E)+(g[1]\otimes
 1+t))$. The fact that $d_{\alpha }^2=0$ implies that $\tilde{g}:=g\otimes
 1+t[-1]:V\to V^\prime$ is a closed morphism of degree zero and hence
 the DG module $C(g\otimes 1)$ with the differential $d_\alpha$ is
 the cone $C(\tilde{g})$ of this morphism.

 Clearly,
 $i^*\tilde{g}=g$ and it remains to prove that $\tilde{g}$ is a
 homotopy equivalence. This in turn is equivalent to the
 acyclicity of the DG algebra $\End(C(\tilde{g}))$. But recall
 that the differential in $\End(C(\tilde{g}))$ is an "$\cR$-deformation" of the
 differential in the  DG algebra
 $\End(C(g))$ which is acyclic, since $g$ is a homotopy
 equivalence. Therefore $\End(C(\tilde{g}))$ is also acyclic. This
 proves the first statement in a). The last statement follows by the
 equivalence of groupoids $\Def ^{\h}_{\cR}\simeq \coDef
 ^{\h}_{\cR}$ (Proposition 4.7).

The proof of b) is similar: exactly in the same way we construct a
closed morphism of degree zero $\tilde{f}:S\to T$ which extends $f$.
Then $\tilde{f}$ is a quasi-isomorphism, because $f$ is such.
\end{proof}

\begin{cor} Fix an artinian DG algebra $\cR$.

a) Let $g:E\to E^\prime $ be a homotopy equivalence as in
Proposition 8.5a). Let $(V,\id )\in \Def ^{\h}_{\cR}(E)$ and
$(V^\prime,\id )\in \Def ^{\h}_{\cR}(E^\prime)$ be objects
corresponding to each other under the equivalence $\Def
^{\h}_{\cR}(E)\simeq \Def ^{\h}_{\cR}(E^\prime)$. Then $i^*V=\bL
i^*V$ if and only if $i^*V^\prime=\bL i^*V^\prime$. Similarly for
the objects of $\coDef _{\cR}^{\h}$ with $i^!$ and $\bR i^!$ instead
of $i^*$ and $\bL i^*$.

b) Let $f:P\to I$ be a quasi-isomorphism as in Proposition 8.5b).
Let $(S,\id)\in \Def _{\cR}^{\h}(P)$ and $(T,\id)\in \Def
_{\cR}^{\h}(I)$ be objects which correspond to each other under the
equivalence $\Def _{\cR}^{\h}(P)\simeq \Def _{\cR}^{\h}(I).$ Then
$i^*S=\bL i^*S$ if and only if $i^*T =\bL i^*T$. Similarly for the
objects of $\coDef _{\cR}^{\h}$ with $i^!$ and $\bR i^!$ instead of
$i^*$ and $\bL i^*$.
\end{cor}

\begin{proof} This follows immediately from Proposition 8.5.
\end{proof}

\begin{prop}  Let $F:\cA \to \cC$ be a DG functor which induces an
equivalence of derived categories $\bL F^*:D(\cA ^{op})\to D(\cC ^{op})$.
(For example, this is the case if $F$ induces a quasi-equivalence
$F^{\pre-tr}:\cA ^{\pre-tr}\to \cC ^{\pre-tr}$ (Corollary 3.15)).

a) Let $P\in \cP (\cA ^{op})$. Then the map of DG algebras $F^*:\End
(P)\to \End (F^*(P))$ is a quasi-isomorphism. Hence the deformation
pseudo-functors $\Def ^{\h}$ and $\coDef ^{\h}$ of $P$ and $F^*(P)$
are equivalent.

b) Let $I\in \cI(\cA ^{op})$. Then the map of DG algebras
$F^!:\End(I)\to \End (F^!(I))$ is a quasi-isomorphism. Hence the
deformation pseudo-functors $\Def ^{\h}$ and $\coDef ^{\h}$ of $I$
and $F^!(I)$ are equivalent.
\end{prop}

\begin{proof} a) By Lemma 3.6 we have $F^*(P)\in \cP (\cC ^{op})$. Hence the assertion follows from
Theorems 3.1 and 8.1.

b) The functor $\bR F^!:D(\cA ^{op})\to D(\cC ^{op})$ is also an
equivalence because of adjunctions $(F_*,\bR F^!),(\bL F^* ,F_*)$.
Also $F^!(I)\in \cI(\cC ^{op})$ (Lemma 3.6). Hence the assertion
follows from Theorems 3.1 and 8.1.
\end{proof}

\section{Direct relation between pseudo-functors $\Def ^{\h}(F)$ and $\Def
^{\h}(\cB)$ ($\coDef ^{\h}(F)$ and $\coDef ^{\h}(\cB)$)}

\subsection{DG functor $\Sigma$} Let $F\in \cA ^{op}\text{-mod}$ and
put $\cB =\End(F)$. Recall the DG functor from Example 3.14
$$\Sigma =\Sigma ^F:\cB ^{op}\text{-mod}\to \cA ^{op}\text{-mod},\quad
\Sigma (M)=M\otimes _{\cB}F.$$ For each artinian DG algebra $\cR$ we
obtain the corresponding DG functor
$$\Sigma _{\cR}:(\cB\otimes \cR) ^{op}\text{-mod}\to \cA _{\cR}^{op}\text{-mod},\quad
\Sigma _{\cR}(M)=M\otimes _{\cB}F.$$

\begin{lemma} The DG functors $\Sigma _{\cR}$ have the following
properties.

a) If a DG $(\cB\otimes \cR) ^{op}$-module $M$ is graded $\cR$-free
(resp. graded $\cR$-cofree), then so is the DG $\cA_{\cR} ^{op}$-module
$\Sigma _{\cR}(M)$.

b) Let $\phi :\cR \to \cQ$ be a homomorphism of artinian DG
algebras. Then there are natural isomorphisms of DG functors
$$\Sigma _{\cQ}\cdot \phi ^*=\phi ^*\cdot \Sigma _{\cR}, \quad
\Sigma _{\cR}\cdot \phi _*=\phi _*\cdot \Sigma _{\cQ}.$$ In
particular,
$$\Sigma \cdot i ^*=i ^*\cdot \Sigma _{\cR}.$$

c) There is a natural isomorphism of DG functors
$$\Sigma _{\cQ}\cdot \phi ^!=\phi ^!\cdot \Sigma _{\cR}$$
on the full DG subcategory of DG $(\cB\otimes \cR)^{op}$-modules $M$
such that $M^{\gr}\simeq M_1^{\gr}\otimes M_2^{\gr}$ for a
$\cB^{op}$-module $M_1$ and an $\cR ^{op}$-module $M_2$. (This subcategory
includes in particular graded $\cR$-cofree modules.) Therefore
$$\Sigma \cdot i ^!=i ^!\cdot \Sigma _{\cR}$$
on this subcategory.

d) For a graded $\cR$-free DG $(\cB\otimes \cR) ^{op}$-module $M$ there
is a functorial isomorphism
$$\Sigma _{\cR}(M\otimes _{\cR} \cR ^*)=\Sigma _{\cR}(M)\otimes _{\cR}\cR ^*$$
\end{lemma}

\begin{proof} The only nontrivial assertion is c). For any DG $(\cB
\otimes \cR)^{op}$-module $M$ there is a natural closed morphism of
degree zero of DG $\cA ^{op}_{\cQ}$-modules
$$\gamma _M:\Hom _{\cR ^{op}}(Q, M)\otimes _{\cB}F\to \Hom _{\cR ^{op}}(Q,
M\otimes _{\cB}F), \quad \gamma(g\otimes
f)(q)=(-1)^{\bar{f}\bar{q}}g(q)\otimes f.$$ Since $\cQ$ is a finite
$\cR ^{op}$-module $\gamma _M$ is an isomorphism if $M^{\gr}\simeq
M_1^{\gr}\otimes M_2^{\gr}$ for a $\cB^{op}$-module $M_1$ and an $\cR
^{op}$-module $M_2$.
\end{proof}

\begin{prop} a) For each artinian DG algebra $\cR$ the DG functor
$\Sigma _{\cR}$ induces functors between groupoids
$$\Def ^{\h}(\Sigma _{\cR}):\Def ^{\h}_{\cR}(\cB)\to \Def ^{\h}_{\cR
}(F),$$
$$\coDef ^{\h}(\Sigma _{\cR}):\coDef ^{\h}_{\cR}(\cB)\to \coDef ^{\h}_{\cR
}(F),$$

b) The collection of DG functors $\{\Sigma _{\cR}\}_{\cR}$ defines
morphisms of pseudo-functors
$$\Def ^{\h}(\Sigma ):\Def ^{\h}(\cB)\to \Def ^{\h}(F),$$
$$\coDef ^{\h}(\Sigma ):\Def ^{\h}(\cB)\to \Def ^{\h}(F).$$

c) The morphism $\Def ^{\h}(\Sigma )$ is compatible with the
equivalence $\theta $ of Proposition 6.1. That is the functorial
diagram
$$\begin{array}{ccc}
\cM \cC (\cB) & = & \cM \cC (\cB)\\
\theta ^{\cB} \downarrow & & \downarrow \theta ^{F}\\
\Def ^{\h}(\cB) & \stackrel{\Def^{\h}(\Sigma)}{\rightarrow} & \Def
^{\h}(F)
\end{array}
$$
is commutative.

d) The morphisms $\Def ^{\h}(\Sigma )$ and $\coDef ^{\h}(\Sigma )$
are compatible with the equivalence $\delta$ of Proposition 4.7.
That is the functorial diagram
$$\begin{array}{ccc}
\Def ^{\h}(\cB) & \stackrel{\Def^{\h}(\Sigma)}{\rightarrow} & \Def
^{\h}(F)\\
\delta ^{\cB} \downarrow & & \downarrow \delta ^{F}\\
\coDef ^{\h}(\cB) & \stackrel{\coDef^{\h}(\Sigma)}{\rightarrow} &
\coDef ^{\h}(F)
\end{array}
$$
is commutative.

e) The morphisms $\Def ^{\h}(\Sigma )$ and $\coDef ^{\h}(\Sigma )$
are equivalences, i.e. for each $\cR$ the functors $\Def
^{\h}(\Sigma _{\cR})$ and $\coDef ^{\h}(\Sigma _{\cR})$ are
equivalences.
\end{prop}

\begin{proof} a) and b) follow from parts a),b),c) of Lemma 9.1; c) is obvious;
d) follows from part d) of Lemma 9.1; e) follows from c) and d).
\end{proof}

\subsection{DG functor $\psi^*$} Let $\psi :\cC \to \cB$ be a
homomorphism of DG algebras. Recall the corresponding DG functor
$$\psi ^*:\cC ^{op}\text{-mod}\to \cB ^{op}\text{-mod},\quad \psi
^*(M)=M\otimes _{\cC}\cB.$$ For each artinian DG algebra $\cR$ we
obtain a similar DG functor
$$\psi ^*_{\cR}:(\cC \otimes \cR) ^{op}
\text{-mod}\to (\cB\otimes \cR) ^{op}\text{-mod},\quad \psi
^*(M)=M\otimes _{\cC}\cB.$$

The next lemma and proposition are complete analogues of Lemma 9.1
and Proposition 9.2.

\begin{lemma} The DG functors $\psi ^* _{\cR}$ have the following
properties.

a) If a DG $(\cC\otimes \cR) ^{op}$-module $M$ is graded $\cR$-free
(resp. graded $\cR$-cofree), then so is the DG $(\cB\otimes \cR)
^{op}$-module $\psi ^*_{\cR}(M)$.

b) Let $\phi :\cR \to \cQ$ be a homomorphism of artinian DG
algebras. Then there are natural isomorphisms of DG functors
$$\psi^* _{\cQ}\cdot \phi ^*=\phi ^*\cdot \psi^* _{\cR}, \quad
\psi^* _{\cR}\cdot \phi _*=\phi _*\cdot \psi ^* _{\cQ}.$$ In
particular,
$$\psi^* \cdot i ^*=i ^*\cdot \psi^* _{\cR}.$$

c) There is a natural isomorphism of DG functors
$$\psi^* _{\cQ}\cdot \phi ^!=\phi ^!\cdot \psi^* _{\cR}$$
on the full DG subcategory of  DG $(\cC\otimes \cR)^{op}$-modules $M$
such that $M^{\gr}\simeq M_1^{\gr}\otimes M_2^{\gr}$ for a
$\cC^{op}$-module $M_1$ and an $\cR ^{op}$-module $M_2$. (This subcategory
includes in particular graded $\cR$-cofree modules.) Therefore
$$\psi^* \cdot i ^!=i ^!\cdot \psi^* _{\cR}$$
on this subcategory.

d) For a graded $\cR$-free DG $(\cC\otimes \cR) ^{op}$-module $M$ there
is a functorial isomorphism
$$\psi^* _{\cR}(M\otimes _{\cR} \cR ^*)=\psi^* _{\cR}(M)\otimes _{\cR}\cR ^*$$
\end{lemma}

\begin{proof} As in Lemma 9.1, the only nontrivial assertion is
c). For any DG $(\cC\otimes \cR)^{op}$-module $M$ there is a natural
closed morphism of degree zero of DG $A_{\cQ}^{op}$-modules
$$\eta_M:\Hom _{\cR ^{op}}(Q, M)\otimes _{\cC}\cB\to \Hom _{\cR ^{op}}(Q,
M\otimes _{\cC}\cB), \quad \gamma(g\otimes
f)(q)=(-1)^{\bar{f}\bar{q}}g(q)\otimes f.$$ Since $\cQ$ is a
finite $\cR ^{op}$-module $\eta_M$ is an isomorphism if
$M^{\gr}\simeq M_1^{\gr}\otimes M_2^{\gr}$ for a $\cB^{op}$-module
$M_1$ and an $\cR ^{op}$-module $M_2$.
\end{proof}

\begin{prop} a) For each artinian DG algebra $\cR$ the DG functor
$\psi^* _{\cR}$ induces functors between groupoids
$$\Def ^{\h}(\psi^* _{\cR}):\Def ^{\h}_{\cR}(\cC)\to \Def ^{\h}_{\cR
}(\cB),$$
$$\coDef ^{\h}(\psi^* _{\cR}):\coDef ^{\h}_{\cR}(\cC)\to \coDef ^{\h}_{\cR
}(\cB),$$

b) The collection of DG functors $\{\psi^* _{\cR}\}_{\cR}$ defines
morphisms
$$\Def ^{\h}(\psi^* ):\Def ^{\h}(\cC)\to \Def ^{\h}(\cB),$$
$$\coDef ^{\h}(\psi^* ):\Def ^{\h}(\cC)\to \Def ^{\h}(\cB).$$

c) The morphism $\Def ^{\h}(\psi^* )$ is compatible with the
equivalence $\theta $ of Proposition 6.1. That is the functorial
diagram
$$\begin{array}{ccc}
\cM \cC (\cC) & \stackrel{\psi^*}{\rightarrow} & \cM \cC (\cB)\\
\theta ^{\cC} \downarrow & & \downarrow \theta ^{\cB}\\
\Def ^{\h}(\cC) & \stackrel{\Def^{\h}(\psi^*)}{\rightarrow} & \Def
^{\h}(\cB)
\end{array}
$$
is commutative.

d) The morphisms $\Def ^{\h}(\psi^* )$ and $\coDef ^{\h}(\psi^* )$
are compatible with the equivalence $\delta$ of Proposition 4.7.
That is the functorial diagram
$$\begin{array}{ccc}
\Def ^{\h}(\cC) & \stackrel{\Def^{\h}(\psi^*)}{\rightarrow} & \Def
^{\h}(\cB)\\
\delta ^{\cC} \downarrow & & \downarrow \delta ^{\cB}\\
\coDef ^{\h}(\cC) & \stackrel{\coDef^{\h}(\psi^*)}{\rightarrow} &
\coDef ^{\h}(\cB)
\end{array}
$$
is commutative.

e) Assume that $\psi$ is a quasi-isomorphism. Then the morphisms
$\Def ^{\h}(\psi^* )$ and $\coDef ^{\h}(\psi^* )$ are equivalences,
i.e. for each $\cR$ the functors $\Def ^{\h}(\psi^* _{\cR})$ and
$\coDef ^{\h}(\psi^* _{\cR})$ are equivalences.

\end{prop}

\begin{proof} a) and b) follow from parts a),b),c) of Lemma 9.3; c)
is obvious; d) follows from part d) of Lemma 9.3; e) follows from
c),d) and Theorem 8.1.
\end{proof}

Later we will be especially interested in the following example.

\begin{lemma}(Keller). a) Assume that the DG algebra $\cB$ satisfies
the following conditions: $H^i(\cB)=0$ for $i<0$, $H^0(\cB)=k$
(resp. $H^0(\cB)=k$). Then there exists a DG subalgebra $\cC\subset
\cB$ with the properties: $\cC^i=0$ for $i<0$, $\cC^0=k$, and the
embedding $\psi:\cC\hookrightarrow \cB$ is a quasi-isomorphism
(resp. the induced map $H^i(\psi):H^i(\cC)\to H^i(\cB)$ is an
isomorphism for $i\geq 0$).
\end{lemma}

\begin{proof} Indeed, put $\cC^0=k$, $\cC^1=K\oplus L$, where
$d(K)=0$ and $K$ projects isomorphically to $H^1(\cB)$, and
$d:L\stackrel{\sim}{\to}d(\cB ^1)\subset \cB ^2$. Then take $\cC
^i=\cB^i$ for $i\geq 2$ and $\cC ^i=0$ for $i<0$.
\end{proof}

\section{The derived deformation and co-deformation pseudo-functors}

\subsection{The pseudo-functor $\Def (E)$} Fix a DG category $\cA$ and an object $E\in \cA^{op}\text{-mod}$. We
are going to define a pseudo-functor $\Def(E)$ from the category
$\dgart$ to the category ${\bf Gpd}$ of groupoids. This
pseudo-functor assigns to a DG algebra $\cR$ the groupoid $\Def
_{\cR}(E)$ of $\cR$-deformations of $E$ in the {\it derived}
category $D(\cA ^{op})$.

\begin{defi} Fix an artinian DG algebra $\cR$. An object of the
groupoid $\Def _{\cR}(E)$ is a pair $(S,\sigma)$, where $S\in D(\cA
_{\cR}^{op})$ and $\sigma$ is an isomorphism (in $D(\cA ^{op})$)
$$\sigma :\bL i^*S\to E.$$
A morphism $f:(S,\sigma)\to (T,\tau)$ between two $\cR$-deformations
of $E$ is an isomorphism (in $D(\cA _{\cR} ^{op})$) $f:S\to T$, such
that
$$\tau \cdot \bL i^*(f)=\sigma.$$
This defines the groupoid $\Def _{\cR}(E)$. A homomorphism of
artinian DG algebras $\phi:\cR \to \cQ$ induces the functor
$$\bL\phi ^*:\Def _{\cR}(E)\to \Def _{\cQ}(E).$$
Thus we obtain a pseudo-functor
$$\Def (E):\dgart \to {\bf Gpd}.$$
\end{defi}

We call $\Def (E)$ the pseudo-functor of derived deformations of
$E$.

\begin{remark} A quasi-isomorphism $\phi :\cR\to \cQ$ of artinian
DG algebras induces an equivalence of groupoids
$$\bL\phi ^*:\Def _{\cR}(E)\to \Def _{\cQ}(E).$$ Indeed,
$\bL\phi ^*:D(\cA _{\cR} ^{op})\to D(\cA _{\cQ} ^{op})$ is an equivalence
of categories (Proposition 3.7) which commutes with the functor $\bL
i^*$.
\end{remark}

\begin{remark} A quasi-isomorphism $\delta:E_1\to E_2$ of
DG $\cA^{op}$-modules induces an equivalence of pseudo-functors
$$\delta _*:\Def (E_1)\to \Def (E_2)$$
by the formula $\delta _*(S,\sigma)=(S,\delta \cdot \sigma)$.
\end{remark}

\begin{prop} Let $F:\cA \to \cA ^\prime$ be a DG functor which
induces a quasi-equivalence  $F^{\pre-tr}:\cA ^{\pre-tr}\to
\cA^{\prime \pre-tr}$ (this happens for example if $F$ is a
quasi-equivalence). Then for any $E\in D(\cA ^{op})$ the deformation
pseudo-functors $\Def (E)$ and $\Def (\bL F^*(E))$ are canonically
equivalent. (Hence also $\Def (F_*(E^\prime))$ and $\Def (E^\prime)$
are equivalent for any $E^\prime \in D(\cA ^{\prime 0})$).
\end{prop}

\begin{proof} For any artinian DG algebra $\cR$ the functor $F$
induces a commutative functorial diagram
$$\begin{array}{ccc}
D(\cA _{\cR} ^{op}) & \stackrel{\bL (F \otimes \id)^*}{\longrightarrow}
&
D(\cA _{\cR}^{\prime 0})\\
\downarrow \bL i^* & & \downarrow \bL i^*\\
D(\cA ^{op}) & \stackrel{\bL F^*}{\longrightarrow} & D(\cA ^{\prime 0})
\end{array}
$$
where $\bL F^*$ and $\bL (F\otimes \id)^*$ are equivalences by
Corollary 3.15.  The horizontal arrows define a functor
$F^*_{\cR}:\Def _{\cR}(E)\to \Def _{\cR}(\bL F^*(E))$. Moreover
these functors are compatible with the functors $\bL \phi ^*:\Def
_{\cR}\to \Def _{\cQ}$ induced by morphisms $\phi :\cR \to \cQ$ of
artinian DG algebras. So we get the morphism $F^*:\Def (E)\to \Def
(\bL F^*(E))$ of pseudo-functors. It is clear that for each $\cR$
the functor $F^*_{\cR}$ is an equivalence. Thus $F^*$ is also such.
\end{proof}

\begin{example} Suppose that $\cA ^\prime$ is a pre-triangulated DG category (so that
the homotopy category $\Ho (\cA ^\prime)$ is triangulated). Let
$F:\cA \hookrightarrow \cA ^\prime$ be an embedding of a full DG
subcategory so that the triangulated category $\Ho (\cA ^\prime)$ is
generated by the collection of objects $F(Ob\cA)$. Then the
assumption of the previous proposition holds.
\end{example}

\begin{remark} In the definition of the pseudo-functor $\Def (E)$ we
could work with the homotopy category of h-projective DG modules
instead of the derived category. Indeed, the functors $i^*$ and
$\phi ^*$ preserve h-projective DG modules.
\end{remark}

\begin{defi} Denote by $\Def _+(E)$, $\Def _-(E)$, $\Def
_0(E)$, $\Def _{\cl}(E)$ the restrictions of the pseudo-functor
$\Def (E)$ to subcategories $\dgart _+$, $\dgart _-$, $\art$,
$\cart$ respectively.
\end{defi}

\subsection{The pseudo-functor $\coDef (E)$}
Now we define the pseudo-functor $\coDef (E)$ of {\it derived
co-deformations} in a similar way replacing everywhere the functors
$(\cdot )^*$ by $(\cdot )^!$.

\begin{defi} Fix an artinian DG algebra $\cR$. An object of the
groupoid $\coDef _{\cR}(E)$ is a pair $(S,\sigma)$, where $S\in
D(\cA _{\cR}^{op})$ and $\sigma$ is an isomorphism (in $D(\cA ^{op})$)
$$\sigma :E\to \bR i^!S.$$
A morphism $f:(S,\sigma)\to (T,\tau)$ between two $\cR$-deformations
of $E$ is an isomorphism (in $D(\cA _{\cR}^{op})$) $f:S\to T$, such
that
$$ \bR i^!(f)\cdot \sigma=\tau.$$
This defines the groupoid $\coDef _{\cR}(E)$. A homomorphism of
artinian DG algebras $\phi:\cR \to \cQ$ induces the functor
$$\bR\phi ^!:\coDef _{\cR}(E)\to \coDef _{\cQ}(E).$$
Thus we obtain a pseudo-functor
$$\coDef (E):\dgart \to {\bf Gpd}.$$
\end{defi}

We call $\coDef (E)$ the functor of derived co-deformations of $E$.

\begin{remark} A quasi-isomorphism $\phi :\cR\to \cQ$ of artinian
DG algebras induces an equivalence of groupoids
$$\bR\phi ^!:\coDef _{\cR}(E)\to \coDef _{\cQ}(E).$$ Indeed,
$\bR\phi ^!:D(\cA _{\cR}^{op})\to D(\cA _{\cQ}^{op})$ is an equivalence of
categories (Proposition 3.7) which commutes with the functor $\bR
i^!$.
\end{remark}

\begin{remark} A quasi-isomorphism $\delta:E_1\to E_2$ of
$\cA$-DG-modules induces an equivalence  of pseudo-functors
$$\delta ^*:\coDef (E_2)\to \coDef (E_1)$$
by the formula $\delta ^*(S,\sigma)=(S,\sigma \cdot \delta)$.
\end{remark}

\begin{prop} Let $F:\cA \to \cA ^\prime$ be a DG functor as in Proposition 10.4 above.
 Consider the induced equivalence of derived
categories $\bR F^!:D(\cA ^{op})\to D(\cA ^{\prime 0})$ (Corollary
3.15). Then for any $E\in D(\cA ^{op})$ the deformation pseudo-functors
$\coDef (E)$ and $\coDef (\bR F^!(E))$ are canonically equivalent.
(Hence also $\coDef (F_*(E^\prime))$ and $\coDef (E^\prime)$ are
equivalent for any $E^\prime \in D(\cA ^{\prime 0})$).
\end{prop}

\begin{proof} For any artinian DG algebra $\cR$ the functor $F$
induces a commutative functorial diagram
$$\begin{array}{ccc}
D(\cA _{\cR}^{op}) & \stackrel{\bR ((F \otimes
\id)^!)}{\longrightarrow} &
D(\cA _{\cR}^{\prime 0})\\
\downarrow \bR i^! & & \downarrow \bR i^!\\
D(\cA ^{op}) & \stackrel{\bR F^!}{\longrightarrow} & D(\cA ^{\prime
0}),
\end{array}
$$
where $\bR (F \otimes \id)^!$ is an equivalence by Corollary 3.15.
The horizontal arrows define a functor $ F^!_{\cR}:\coDef
_{\cR}(E)\to \coDef _{\cR}(\bR F^!(E))$. Moreover these functors are
compatible with the functors $\bR \phi ^!:\coDef _{\cR}\to \coDef
_{\cQ}$ induced by morphisms $\phi :\cR \to \cQ$ of artinian DG
algebras. So we get the morphism $F^!:\coDef (E)\to \coDef (\bR
F^!(E))$. It is clear that for each $\cR$ the functor $F^!_{\cR}$ is
an equivalence. Thus $F^!$ is also such.
\end{proof}

\begin{example} Let $F:\cA ^\prime \to \cA $ be as in Example 10.5 above.
Then the assumption of the previous proposition holds.
\end{example}

\begin{remark} In the definition of the pseudo-functor $\coDef (E)$ we
could work with the homotopy category of h-injective DG modules
instead of the derived category. Indeed, the functors $i^!$ and
$\phi ^!$ preserve h-injective DG modules.
\end{remark}

\begin{defi} Denote by $\coDef _+(E)$, $\coDef _-(E)$, $\coDef
_0(E)$, $\coDef _{\cl}(E)$ the restrictions of the pseudo-functor
$\coDef (E)$ to subcategories $\dgart _+$, $\dgart _-$, $\art$,
$\cart$ respectively.
\end{defi}

\begin{remark} The pseudo-functors $\Def (E)$ and $\coDef (E)$ are not always equivalent
(unlike their homotopy counterparts $\Def ^{\h}(E)$ and $\coDef
^{\h}(E)$). In fact we expect that pseudo-functors $\Def$ and
$\coDef $ are the "right ones" only in case they can be expressed in
terms of the pseudo-functors $\Def ^{\h}$ and $\coDef ^{\h}$
respectively. (See the next section).
\end{remark}

\section{Relation between pseudo-functors $\Def$ and $\Def ^h$ (resp.
$\coDef $ and $\coDef ^h$)}

The ideal scheme that should relate these deformation
pseudo-functors is the following. Let $\cA$ be a DG category, $E\in
\cA ^{op}\text{-mod}$. Choose quasi-isomorphisms $P\to E$ and $E\to I$,
where $P\in \P(\cA ^{op})$ and $I\in \cI(\cA ^{op})$. Then there should
exist natural equivalences
$$\Def(E)\simeq \Def ^h (P),\quad\quad \coDef(E)\simeq \coDef
^h(I).$$ Unfortunately, this does not always work.

\begin{example} Let $\cA$ be just a graded algebra $A=k[t]$, i.e.
$\cA$ contains a single object with the endomorphism algebra
$k[t]$, $\deg (t)=1$ (the differential is zero). Take the artinian
DG algebra $\cR$ to be $\cR=k[\epsilon]/(\epsilon ^2)$,
$\deg(\epsilon)=0$. Let $E=A$ and consider a DG $\cA
_{\cR}^{op}$-module $M=E\otimes \cR$ with the differential $d_M$
which is the multiplication by $t\otimes \epsilon$. Clearly, $M$
defines an object in $\Def ^h_{\cR}(E)$ which is not isomorphic to
the trivial deformation. However, one can check (Proposition
11.18) that $\bL i^*M$ is not quasi-isomorphic to $E$ (although
$i^*M=E$), thus $M$ does not define an object in $\Def _{\cR}(E)$.
This fact and the next proposition show that the groupoid $\Def
_{\cR}(E)$ is connected (contains only the trivial deformation),
so it is not the "right" one.
\end{example}

\begin{prop} Assume that $\Ext ^{-1}(E,E)=0$.

1) Fix a quasi-isomorphism $P\to E$, $P\in \cP(\cA ^{op})$. Let $\cR $
be an artinian DG algebra and $(S, \id)\in \Def ^{\h}_{\cR}(P)$. The
following conditions are equivalent:

a) $S\in \P(\cA _{\cR}^{op})$,

b) $i^*S=\bL i^*S$,

c) $(S,\id)$ defines an object in the groupoid $\Def _{\cR}(E)$.

The pseudo-functor $\Def (E)$ is equivalent to the full
pseudo-subfunctor of $\Def ^{\h}(P)$ consisting of objects $(S,\id)
\in  \Def ^{\h}(P)$, where $S$ satisfies a) (or b)) above.

2) Fix a quasi-isomorphism  $E\to I$ with $I\in \cI(\cA ^{op})$. Let
$\cR $  be an artinian DG algebra and $(T, \id)\in \coDef ^{\h}
_{\cR}(I)$. The following conditions are equivalent:

a') $T\in \cI(\cA _{\cR}^{op})$,

b') $i^!T=\bR i^!T$,

c') $(T,\id)$ defines an object in the groupoid $\coDef _{\cR}(E)$.

The pseudo-functor $\coDef (E)$ is equivalent to the full
pseudo-subfunctor of $\coDef ^{\h}(I)$ consisting of objects
$(T,\id) \in  \coDef ^{\h}(I)$, where $T$ satisfies a') (or b'))
above.
  \end{prop}

\begin{proof} 1) It is clear that a) implies b) and b) implies c).
We will prove that c) implies a). We may and will replace the
pseudo-functor $\Def (E)$ by an equivalent pseudo-functor $\Def (P)$
(Remark 10.3).

Since $(S, \id)$ defines an object in $\Def _{\cR}(P)$ there exists
a quasi-isomorphism $g:\tilde{S}\to S$ where $\tilde{S}$ has
property (P) (hence $\tilde{S}\in \P(\cA _{\cR}^{op})$), such that
$i^*g: i^*\tilde{S}\to i^*S=P$ is also a quasi-isomorphism. Denote
$Z=i^*\tilde{S}$. Then $Z\in \P(\cA ^{op})$ and hence $i^*g$ is a
homotopy equivalence. Since both $\tilde{S}$ and $S$ are graded
$\cR$-free, the map $g$ is also a homotopy equivalence (Proposition
3.12d)). Thus $S\in \P(\cA _{\cR}^{op})$.

Let us prove the last assertion in 1).

 Fix an object $(\overline{S} ,\tau) \in \Def _{\cR}(P)$. Replacing $(\overline{S},
\tau)$ by an isomorphic object we may and will assume that
$\overline{S}$ satisfies property (P). In particular,
$\overline{S}\in \cP(\cA _{\cR}^{op})$ and $\overline{S}$ is graded
$\cR$-free. This implies that $(\overline{S},\id)\in \Def^{\h}
_{\cR}(W)$ where $W=i^*\overline{S}$. We have $W\in \P(\cA ^{op})$. The
quasi-isomorphism $\tau :W\to P$ is therefore a homotopy
equivalence. By Corollary 8.4a) and Proposition 8.5a) there exists
an object $(S^\prime,\id)\in \Def ^{\h}_{\cR}(P)$ and a homotopy
equivalence $\tau ^\prime :\overline{S}\to S^\prime$ such that
$i^*(\tau ^\prime)=\tau$. This shows that $(\overline{S}, \tau)$ is
isomorphic (in $\Def _{\cR}(P)$) to an object $(S^\prime ,\id )\in
\Def ^{\h}_{\cR}(P)$, where $S^\prime \in \P(\cA _{\cR}^{op})$.

Let $(S ,\id ),(S^\prime, \id) \in \Def _{\cR}^{\h}(P)$ be two
objects such that $S,S^\prime \in \cP(\cA _{\cR}^{op})$.
 Consider the
obvious map
$$\delta :\Hom _{\Def ^{\h}_{\cR}(P)}((S,\id ),(S^\prime,\id ))\to
\Hom _{\Def _{\cR}(P)}((S,\id),(S^\prime,\id)).$$ It suffices to
show that $\delta$ is bijective.

Let $f:(S,\id ) \to (S^\prime, \id )$ be an isomorphism in $\Def
_{\cR}(P)$. Since $S,S^\prime\in \cP(\cA _{\cR}^{op})$ and $P \in
\cP(\cA ^{op})$ this isomorphism $f$ is a homotopy equivalence $f:S\to
S^\prime$ such that $i^*f$ is homotopic to $\id _P$. Let $h:i^*f\to
\id$ be a homotopy. Since $S$, $S^\prime$ are graded $\cR$-free the
map $i^*:\Hom (S, S^\prime )\to \Hom (P,P)$ is surjective
(Proposition 3.12a)). Choose a lift $\tilde{h}:S\to S^\prime[1]$ of
$h$ and replace $f$ by $\tilde{f}=f-d\tilde{h}$. Then
$i^*\tilde{f}=id$. Since $S$ and $S^\prime$ are graded $\cR$-free
$\tilde{f}$ is an isomorphism (Proposition 3.12d)). This shows that
$\delta$ is surjective.

Let $g_1,g_2:S\to S^\prime$ be two isomorphisms (in $\cA
_{\cR}^{op}\text{-mod}$) such that $i^*g_1=i^*g_2=\id _P$. That is
$g_1,g_2$ represent morphisms in $\Def _{\cR}^{\h}(P)$. Assume that
$\delta (g_1)=\delta (g_2)$, i.e. there exists a homotopy $s:g_1\to
g_2$. Then $d(i^*s)=i^*(ds)=0$. Since by our assumption $H^{-1}\Hom
(P,P)=0$ there exists $t\in \Hom ^{-2}(P,P)$ with $dt=i^*s$. Choose
a lift $\tilde{t}\in \Hom ^{-2}(S,S^\prime)$ of $t$. Then
$\tilde{s}:=s-d\tilde{t}$ is an allowable homotopy between $g_1$ and
$g_2$. This proves that $\delta $ is injective and finishes the
proof of 1).

The proof of 2) is very similar, but we present it for completeness.
Again it is clear that a') implies b') and b') implies c'). We will
prove that c') implies a') We may and will replace the functor
$\coDef (E)$ by an equivalent functor $\coDef (I)$ (Remark 10.10).

Since $(T,\id)$ defines an object in $\coDef _{\cR}(I)$, there
exists a quasi-isomorphism $g:T\to \tilde{T}$ where $\tilde{T}$ has
property (I) (hence $\tilde{T} \in \cI(\cA _{\cR}^{op})$), such that
$i^!g:I=i^!T\to i^!\tilde{T}$ is also a quasi-isomorphism. Denote
$K=i^!\tilde{T}$. Then $K\in \cI(\cA ^{op})$ and hence $i^!g$ is a
homotopy equivalence. Since both $T$ and $\tilde{T}$ are graded
$\cR$-cofree, the map $g$ is also a homotopy equivalence
(Proposition 3.12d)). Thus $T\in \cI(\cA _{\cR} ^{op})$.

Let us prove the last assertion in 2).

 Fix an object $(\overline{T} ,\tau) \in \coDef _{\cR}(I)$.
 Replacing $(\overline{T}, \tau)$ by an isomorphic object we may and will
assume that $\overline{T}$ satisfies property (I). In particular,
$\overline{T}\in \cI(\cA _{\cR} ^{op})$ and $\overline{T}$ is graded
$\cR$-cofree. This implies that $(\overline{T},\id)\in \coDef^{\h}
_{\cR}(L)$ where $L=i^!\overline{T}$. We have $L\in \cI (\cA ^{op})$
and hence the quasi-isomorphism $\tau :I\to L$ is a homotopy
equivalence. By Corollary 8.4a) and Proposition 8.5a) there exist an
object $(T^\prime,\id)\in \coDef ^{\h}_{\cR}(I)$ and a homotopy
equivalence $\tau ^\prime : T^\prime\to \overline{T}$ such that
$i^!\tau ^\prime =\tau$. In particular, $T^\prime\in \cI(\cA _{\cR}
^{op})$. This shows that $(\overline{T},\tau)$ is isomorphic (in
$\coDef _{\cR}(I)$) to an object $(T^\prime ,\id) \in \coDef
_{\cR}^{\h}(I)$ where $T^\prime \in \cI(\cA _{\cR}^{op})$.

Let $(T ,\id ),(T^\prime, \id) \in \coDef _{\cR}^{\h}(I)$ be two
objects such that $T,T^\prime \in \cI(\cA _{\cR}^{op})$.
 Consider the
obvious map
$$\delta :\Hom _{\coDef ^{\h}_{\cR}(I)}((T,\id ),(T^\prime,\id ))\to
\Hom _{\coDef _{\cR}(I)}((T,\id),(T^\prime,\id)).$$ It suffices to
show that $\delta$ is bijective.

Let $f:(T,\id ) \to (T^\prime, \id )$ be an isomorphism in $\coDef
_{\cR}(I)$. Since $T,T^\prime\in \cI(\cA _{\cR}^{op})$ and $I \in
\cI(\cA ^{op})$ this isomorphism $f$ is a homotopy equivalence $f:T\to
T^\prime$ such that $i^!f$ is homotopic to $\id _I$. Let $h:i^!f\to
\id$ be a homotopy.  Since $T$, $T^\prime$ are graded $\cR$-cofree
the map $i^!:\Hom (T, T^\prime )\to \Hom (I,I)$ is surjective
(Proposition 3.12a)). Choose a lift $\tilde{h}:T\to T^\prime[1]$ of
$h$ and replace $f$ by $\tilde{f}=f-d\tilde{h}$. Then
$i^!\tilde{f}=id$. Since $T$ and $T^\prime$ are graded $\cR$-cofree
$\tilde{f}$ is an isomorphism (Proposition 3.12d)). This shows that
$\delta$ is surjective.

Let $g_1,g_2:T\to T^\prime$ be two isomorphisms (in $\cA
_{\cR}^{op}\text{-mod}$) such that $i^!g_1=i^!g_2=\id _I$. That is
$g_1,g_2$ represent morphisms in $\coDef _{\cR}^{\h}(I)$. Assume
that $\delta (g_1)=\delta (g_2)$, i.e. there exists a homotopy
$s:g_1\to g_2$. Then $d(i^!s)=i^!(ds)=0$. Since by our assumption
$H^{-1}\Hom (I,I)=0$ there exists $t\in \Hom ^{-2}(I,I)$ with
$dt=i^!s$. Choose a lift $\tilde{t}\in \Hom ^{-2}(T,T^\prime)$ of
$t$. Then $\tilde{s}:=s-d\tilde{t}$ is an allowable homotopy between
$g_1$ and $g_2$. This proves that $\delta $ is injective.
\end{proof}

\begin{remark} In the situation of Proposition 11.2
using Corollary 8.4b)  also obtain full and faithful morphisms of
pseudo-functors $ \Def (E)$, $\coDef (E)$ to each of the equivalent
pseudo-functors $\Def ^{\h}(P)$, $\coDef ^{\h}(P)$, $\Def ^{\h}(I)$,
$\coDef ^{\h}(I)$.
\end{remark}

\begin{cor} Assume that $\Ext ^{-1}(E,E)=0$. Let $F\in \cA
^{op}\text{-mod}$ be an h-projective or an h-injective quasi-isomorphic
to $E$.

a) The pseudo-functor $\Def (E)$ ($\simeq \Def (F)$) is equivalent
to the full pseudo-subfunctor of $\Def ^{\h}(F)$ which consists of
objects $(S,\id )$ such that $i^*S=\bL i^*S$.

b) The pseudo-functor $\coDef (E)$ ($\simeq \coDef (F))$ is
equivalent to the full pseudo-subfunctor of $\coDef ^{\h}(F)$ which
consists of objects $(T,\id)$ such that $i^!T=\bR i^!T$.
\end{cor}

\begin{proof} a). In case $F$ is h-projective this is Proposition
11.2 1). Assume that $F$ is h-injective. Choose a
quasi-isomorphism $P\to F$ where $P$ is h-projective. Again by
Proposition 11.2 1) the assertion holds for $P$ instead of $F$.
But then it also holds for $F$ by Corollary 8.6 b).

b). In case $F$ is h-injective this is Proposition 11.2 2). Assume
that $F$ is h-projective. Choose a quasi-isomorphism $F\to I$ where
$I$ is h-injective. Then again by Proposition 11.2 2) the assertion
holds for $I$ instead of $F$. But then it also holds for $F$ by
Corollary 8.6 b).
\end{proof}

The next theorem provides an example when the pseudo-functors $\Def
_-$ and $\Def ^{\h}_-$ (resp. $\coDef _- $ and $\coDef ^{\h}_-$) are
equivalent.

\begin{defi} An object $M\in \cA ^{op}\text{-mod}$ is called bounded above (resp. below) if
there exists $i$ such that $M(A)^j=0$ for all $A\in \cA$ and all
$j\geq i$ (resp. $j\leq i$).
\end{defi}

\begin{theo} Assume  that $\Ext ^{-1}(E,E)=0$.

a) Suppose that there exists an h-projective or an h-injective $P\in
\cA ^{op}\text{-mod}$ which is  bounded above  and  quasi-isomorphic to
$E$. Then the pseudo-functors $\Def _-(E)$ and $\Def _-^{\h}(P)$ are
equivalent.

 b) Suppose that there exists an h-projective or an h-injective $I\in
\cA ^{op}\text{-mod}$ which is  bounded below  and  quasi-isomorphic to
$E$. Then the pseudo-functors $\coDef _-(E)$ and $\coDef _-^{\h}(I)$
are equivalent.
\end{theo}

\begin{proof} Fix $\cR \in \dgart _-$.
In both cases it suffices to show that the embedding of groupoids
$\Def _{\cR} (E)\simeq \Def _{\cR}(P) \subset \Def _{\cR}^{\h}(P)$
(resp. $\coDef _{\cR} (E)\simeq \coDef _{\cR}(I) \subset \coDef
_{\cR}^{\h}(I)$) in Corollary 11.4 is essentially surjective.

a)   It suffices to prove the following lemma.

\begin{lemma} Let $M\in \cA ^{op}\text{-mod}$ be bounded above and $(S,\id)\in \Def _{\cR}^{\h}(M)$.
The DG $\cA _{\cR}^{op}$-module $S$ is acyclic for the functor $i^*$,
i.e. $\bL i^*S=i^*S$.
\end{lemma}

Indeed, in case $M=P$ the lemma implies that $S$ defines an object
in $\Def _{\cR}(P)$ (Corollary 11.4 a)).

\begin{proof} Choose a quasi-isomorphism $f:Q\to S$ where $Q\in
\cP(\cA _{\cR} ^{op})$. We need to prove that $i^*f$ is a
quasi-isomorphism. It suffices to prove that $\pi _!i^*f$ is a
quasi-isomorphism (Example 3.13). Recall that $\pi _!i^*=i^*\pi
_!$. Thus it suffices to prove that $\pi _!f$ is a homotopy
equivalence. Clearly $\pi _!f$ is a quasi-isomorphism. The DG
$\cR^{op}$-module $\pi _!Q$ is h-projective (Example 3.13). We claim
that the DG $\cR ^{op}$-module $\pi _!S $ is also h-projective. Since
the direct sum of h-projective DG modules is again h-projective,
it suffices to prove that for each object $A\in \cA$ the DG
$\cR^{op}$-module $S(A)$ is h-projective. Take some object $A\in
\cA$. We have that $S(A)$ is bounded above and since $\cR \in
\dgart _-$ this DG $\cR ^{op}$-module has an increasing filtration
with subquotients being free DG $\cR ^{op}$-modules. Thus $S(A)$
satisfies property (P) and hence is h-projective.  It follows that
the quasi-isomorphism $\pi _!f:\pi_! Q\to \pi _!S$ is a homotopy
equivalence. Hence $i^*\pi _!f=\pi _!i^*f$ is also such.
\end{proof}

b) The following lemma implies (by Corollary 11.4 b)) that an object
in $\coDef _{\cR}^{\h}(I)$ is also an object in $\coDef _{\cR}(I)$,
which proves the theorem.
\end{proof}

\begin{lemma} Let $T\in \cA _{\cR}^{op}\text{-mod}$ be graded cofree
and bounded below. Then  $T$ is acyclic for the functor $i^!$, i.e.
$\bR i ^!T=i ^!T$.
\end{lemma}

\begin{proof} Denote $N=i^!T\in \cA ^{op}\text{-mod}$. Choose a quasi-isomorphism $g:T\to J$ where $J\in \cI(\cA
_{\cR} ^{op})$. We need to prove that $i^!g$ is a quasi-isomorphism. It
suffices to show that $\pi _* i^!g$ is a quasi-isomorphism. Recall
that $\pi _*i^!=i^!\pi _*$. Thus it suffices to prove that $\pi _*g$
is a homotopy equivalence. Clearly it is a quasi-isomorphism.

Recall that the DG $\cR^{op}$-module $\pi _*J$ is h-injective
(Example 3.13) We claim that $\pi _*T$ is also such. Since the
direct product of h-injective DG modules is again h-injective, it
suffices to prove that for each object $A\in \cA$ the DG
$\cR^{op}$-module $T(A)$ is h-injective. Take some object $A\in \cA$.
Since $\cR \in \dgart _-$ the DG $\cR ^{op}$-module $T(A)$ has a
decreasing filtration
$$G^0\supset G^1\supset G^2\supset
...,$$ with
$$\gr T(A)=\oplus _{j}(T(A))^j \otimes \cR ^*.$$
A direct sum of shifted copies of the DG $\cR ^{op}$-module $\cR ^*$
is h-injective (Lemma 3.18). Thus each $(T(A))^j \otimes \cR ^*$
is h-injective and hence each quotient $T(A)/G^j$ is h-injective.
Also
$$T(A)=\lim_{\leftarrow}T(A)/G^j.$$
Therefore $T(A)$ is h-injective by Remark 3.5.

It follows that $\pi _*g$ is a homotopy equivalence, hence also
$i^!\pi _*g$ is such.
\end{proof}

The last theorem allows us to compare the functors $\Def _-$ and
$\coDef _-$ in some important special cases. Namely we have the
following corollary.

\begin{cor} Assume that

a) $\Ext^{-1}(E,E)=0$;

b) there exists an h-projective or an h-injective $P\in \cA
^{op}\text{-mod}$ which is  bounded above
 and  quasi-isomorphic to $E$;

c)  there exists an h-projective or an h-injective $I\in \cA
^{op}\text{-mod}$ which is  bounded below
 and  quasi-isomorphic to $E$;

Then the pseudo-functors $\Def _-(E)$  and $\coDef _-(E)$ are
equivalent.
\end{cor}

\begin{proof} We have a quasi-isomorphism $P\to I$. Hence by Proposition 8.3 the DG algebras
$\End (P)$ and $\End (I)$ are quasi-isomoprhic. Therefore, in
particular, the pseudo-functors $\Def _-^{\h}(P)$ and $\coDef
_-^{\h}(I)$ are equivalent (Corollary 8.4b)). It remains to apply
the last theorem.
\end{proof}

In practice in order to find the required bounded resolutions one
might need to pass to a "smaller" DG category. So it is useful to
have the following stronger corollary.

\begin{cor} Let $F:\cA \to \cA ^\prime$ be a DG functor which induces a
quasi-equivalence $F^{\pre-tr}:\cA ^{\pre-tr}\to \cA ^{\prime
\pre-tr}$. Consider the corresponding equivalence $F_*:D(\cA
^{\prime 0})\to D(\cA ^{op})$ (Corollary 3.15). Let $E\in \cA ^{\prime
0}\text{-mod}$ be such that

a) $\Ext^{-1}(E,E)=0$;

b) there exists an h-projective or an h-injective $P\in \cA
^{op}\text{-mod}$ which is  bounded above
 and  quasi-isomorphic to $F_*(E)$;

c) there exists an h-projective or an h-injective $P\in \cA
^{op}\text{-mod}$ which is  bounded below
 and  quasi-isomorphic to $F_*(E)$;

 Then the pseudo-functors $\Def _-(E)$ and $\coDef _-(E)$ are
equivalent.
\end{cor}

\begin{proof} By the above corollary the pseudo-functors
$\Def _-(F_*(E))$  and $\coDef _-(F_*(E))$ are equivalent. By
Proposition 10.4 the pseudo-functors $\Def _-(E)$ and $\Def
_-(F_*(E))$ are equivalent. Since the functor $\bR F^!:D(\cA ^{op})\to
D(\cA ^{\prime 0})$ is also an equivalence, we conclude that the
pseudo-functors $\coDef _-(E)$ and $\coDef _-(F_*(E))$ are
equivalent by Proposition 10.11.
\end{proof}

\begin{example} If in the above corollary the DG category $\cA
^\prime$ is pre-triangulated, then one can take for $\cA$ a full DG
subcategory of $\cA ^\prime$ such that $\Ho (\cA ^\prime)$ is
generated as a triangulated category by the subcategory $\Ho (\cA)$.
One can often choose $\cA$ to have one object.
\end{example}

\begin{example} Let $\cC$ be a bounded DG algebra, i.e. $\cC ^i=0$
for $|i|>>0$ and also $H^{-1}(\cC)=0$. Then by Theorem 11.6 and
Proposition 4.7
$$\coDef _-(\cC)\simeq \coDef ^{\h}_-(\cC)\simeq \Def ^{\h}_-(\cC)\simeq \Def _-(\cC).$$
\end{example}

The following theorem makes the equivalence of Corollary 11.9 more
explicit. Let us first introduce some notation.

For an artinian DG algebra $\cR$ consider the DG functors
$$\eta _{\cR}, \epsilon _{\cR}: \cA _{\cR}^{op}\text{-mod} \to \cA
_{\cR}^{op}\text{-mod}$$ defined by
$$\epsilon _{\cR} (M)=M\otimes _{\cR}\cR ^*, \quad \eta
_{\cR}(N)=\Hom _{\cR ^{op}}(\cR ^* ,N).$$ They induce the corresponding
functors
$$\bR \eta _{\cR}, \bL \epsilon _{\cR}:D(\cA ^{op}_{\cR})\to D(\cA
^{op}_{\cR}).$$

\begin{theo} Let $E \in \cA ^{op}\text{-mod}$ satisfy the assumptions
a), b), c) of Corollary 11.9. Fix $\cR \in \dgart _-$. Then the
following holds.

1) Let $F\in \cA ^{op}\text{-mod}$ be h-projective or h-injective
quasi-isomorphic to $E$.

a) For any $(S, \sigma )\in \Def ^{\h}_{\cR}(F)$ we have $i^*S=\bL
i^*S$.

b) For any $(T,\tau )\in \coDef ^{\h}_{\cR}(F)$ we have $i^!T=\bR
i^!T$.

2) There are natural equivalences of pseudo-functors $\Def
^{\h}_-(F)\simeq \Def _-(E)$, $\coDef _-^{\h}(F)\simeq \coDef
_-(E)$.

3) The functors $\bL \epsilon _{\cR}$ and $\bR \eta _{\cR}$ induce
mutually inverse equivalences
$$\bL \epsilon _{\cR}:\Def _{\cR}(E)\to \coDef _{\cR}(E),$$
$$\bR \eta _{\cR}:\coDef _{\cR}(E)\to \Def _{\cR}(E).$$
\end{theo}

\begin{proof} 1a). We may and will assume that $\sigma =\id$.

 Choose a bounded above h-projective or h-injective $P\in \cA ^{op}\text{-mod}$,
 which is quasi-isomorphic to $E$. Then there
exists a quasi-isomorphism $P\to F$ (or $F\to P$). The
pseudo-functors $\Def _-^{\h}(P)$ and $\Def _-^{\h}(F)$ are
equivalent  by Corollary 8.4 (a) or b)). By Theorem 11.6 a) $\Def
_-^{\h}(P)\simeq \Def _-(P)$. Hence by Corollary 11.4 a) for each
$(S^\prime ,\id)\in \Def _{\cR}(P)$ we have $i^*S^\prime =\bL
i^*S^\prime $. Now Corollary 8.6 (a) or b)) implies that $i^*S=\bL
i^*S$. This proves 1a).

1b). We may and will assume that $\tau =\id$.

The proof is similar to that of 1a). Namely, choose a bounded below
h-projective or h-injective $I\in \cA ^{op}\text{-mod}$
quasi-isomorphic to $E$. Then there exists a quasi-isomorphism $F\to
I$ (or $I\to F$).  The pseudo-functors $\coDef _-^{\h}(I)$ and
$\coDef _-^{\h}(F)$ are equivalent and by Corollary 8.4 (a) or b)).
By Theorem 11.6 a) $\coDef _-^{\h}(I)\simeq \coDef _-(I)$. Hence by
Corollary 11.4 b) for each $(T^\prime ,\id)\in \coDef ^{\h}(I)$ we
have $i^!T^\prime =\bR i^!T^\prime$. Now Corollary 8.6 (a) or b))
implies that $i^!T=\bR i^!T$.

2) This follows from 1), Corollary 11.4 a), b).

3) This follows from 2) and the fact that $\epsilon _{\cR}$ and
$\eta _{\cR}$ induce inverse equivalences between $\Def
_{\cR}^{\h}(F)$ and $\coDef _{\cR}^{\h}(F)$ (Proposition 4.7).
\end{proof}

\begin{prop} Let DG algebras $\cB$ and $\cC$ be quasi-isomorphic and
$H^{-1}(\cB)=0$ ($=H^{-1}(\cC)$).  Suppose that the pseudo-functors
$\Def (\cB)$ and $\Def ^{\h}(\cB)$ (resp. $\coDef (\cB)$ and $\coDef
^{\h}(\cB)$) are equivalent. Then the same is true for $\cC$.

Similar results hold for the pseudo-functors $\Def _-, \Def ^{\h}
_-, \coDef _-, ...$.
\end{prop}

\begin{proof}
We may and will assume that there exists a morphism of DG algebras
$\psi :\cB \to \cC$ which is a quasi-isomorphism.

By Proposition 8.6 a) the pseudo-functors $\Def ^{\h}(\cB)$ and
$\Def ^{\h}(\cC)$ are equivalent.

By Proposition 10.4 the pseudo-functors $\Def (\cB)$ are $\Def
(\cC)$ are equivalent.

By Proposition 11.2 a) $\Def (\cB)$ (resp. $\Def (\cC)$) is a full
pseudo-subfunctor  of $\Def ^{\h}(\cB)$ (resp. $\Def ^{\h}(\cC)$).

Thus is $\Def (\cB)\simeq \Def ^{\h}(\cB)$, then also $\Def
(\cC)\simeq \Def ^{\h}(\cC)$.

The proof for $\coDef $ and $\coDef ^{\h}$ is similar using
Proposition 8.6 a), Proposition 10.11 and Proposition 11.2 b).
\end{proof}

\begin{cor} Let $\cB$ be a DG algebra such that $H^{-1}(\cB)=0$.
Assume that $\cB$ is quasi-isomorphic to a DG algebra $\cC$ such
that $\cC$ is bounded above (resp. bounded below). Then the
pseudo-functors $\Def _-(\cB)$ and $\Def ^{\h}_-(\cB)$ are
equivalent (resp. $\coDef _-(\cB)$ and $\coDef ^{\h}_-(\cB)$ are
equivalent).
\end{cor}

\begin{proof} By Theorem 11.6 a) we have that $\Def _-(\cC)$ and $\Def _-^{\h}(\cC)$ are equivalent (resp.
$\coDef _-(\cC)$ and $\coDef _-^{\h}(\cC)$ are equivalent). It
remains to apply Proposition 11.14.
\end{proof}

\subsection{Relation between pseudo-functors $\Def _-(E)$, $\coDef
_-(E)$ and $\Def _-(\cC)$, $\coDef _-(\cC)$}

The next proposition  follows immediately from our previous results.

\begin{prop} Let $\cA$ be a DG category and $E\in \cA ^{op}\text{-mod}$. Assume that

a) $\Ext ^{-1}(E,E)=0$;

b) there exists a bounded above (resp. bounded below) h-projective
or h-injective $F\in \cA ^{op}\text{-mod}$ which is quasi-isomorphic to
$E$;

c) there exists a bounded above (resp. bounded below) DG algebra
$\cC$ which is quasi-isomorphic to $\End (F)$.

Then the pseudo-functors $\Def _-(E)$ and $\Def _-(\cC)$ (resp.
 $\coDef _-(E)$ and $\coDef _- (\cC)$) are equivalent.
\end{prop}

\begin{proof} Assume that $F$ and $\cC$ are bounded above. Then
$\Def _-(E)\simeq \Def ^{\h}_-(F)$ and $\Def _-(\cC)\simeq \Def
^{\h}_-(\cC)$ by Theorem 11.6 a). Also $\Def ^{\h}_-(F)\simeq \Def
^{\h}_-(\cC)$ by Proposition 6.1 and Theorem 8.1.

Assume that $F$ and $\cC$ are bounded below. Then $\coDef
_-(E)\simeq \coDef ^{\h}_-(F)$ and $\coDef _-(\cC)\simeq \coDef
^{\h}_-(\cC)$ by Theorem 11.6 b). Also $\coDef ^{\h}_-(F)\simeq
\coDef ^{\h}_-(\cC)$ by Proposition 6.1 and Theorem 8.1.
\end{proof}

\begin{remark}  The equivalences of pseudo-functors
$\Def ^{\h}_-(\cC)\simeq \Def ^{\h}_-(F)$, $\coDef
^{\h}_-(\cC)\simeq \coDef ^{\h}_-(F)$ in the proof of last
proposition can be made explicit. Put $\cB =\End(F)$. Assume, for
example, that $\psi :\cC \to \cB$ is a homomorphism of DG algebras
which is a quasi-isomorphism. Then the composition of DG functors
(Propositions 9.2, 9.4)
$$\Sigma ^F\cdot \psi ^*:\cC ^{op}\text{-mod}\to \cA ^{op}\text{-mod}$$
induces equivalences of pseudo-functors
$$\Def ^{\h} (\Sigma ^F \cdot \psi ^*): \Def ^{\h}(\cC)\simeq \Def
^{\h}(F)$$
$$\coDef ^{\h} (\Sigma ^F \cdot \psi ^*):\coDef ^{\h}(\cC)\simeq \coDef
^{\h}(F)$$  by Propositions 9.2e) and 9.4f).
\end{remark}

\subsection{Pseudo-functors $\Def(E)$, $\coDef(E)$ are not
determined by the DG algebra $\bR \Hom (E,E)$}

One might expect that the derived deformation and co-deformation
pseudo-functors $\Def _-(E)$, $\coDef _-(E)$ depend only on the
(quasi-isomorphism class of the) DG algebra $\bR \Hom (E,E)$. This
would be an analogue of Theorem 8.1 for the derived deformation
theory. Unfortunately this is not true as is shown in the next
proposition (even for the "classical" pseudo-functors $\Def _{\cl}$,
$\coDef _{\cl}$). This is why all our comparison results for the
pseudo-functors $\Def _-$ and $\coDef _-$ such as
 Theorems 11.6, 11.13, Corollaries 11.9, 11.15, Proposition 11.16 need some boundedness
assumptions.

Consider the DG algebra $A=k[x]$ with the zero differential and
$\deg (x)=1$. Let $\cA$ be the DG category with one object whose
endomorphism DG algebra is $A$. Then $\cA ^{op}\text{-mod}$ is the DG
category of DG modules over the DG algebra $A^{op}=A$. Denote by abuse
of notation the unique object of $\cA$ also by $A$ and consider the
DG $\cA ^{op}$-modules $P=h^A$ and $I=h_A^*$. The first one is
h-projective and bounded below while the second one is h-injective
and bounded above (they are the graded dual of each other). Note
that the DG algebras $\End (P)$ and $\End (I)$ are isomorphic:
$$\End (P)=A, \quad \End (I)=A^{**}=A.$$

Let $\cR =k[\epsilon]/(\epsilon ^2)$ be the (commutative) artinian
DG algebra with the zero differential and $\deg(\epsilon)=0$.

\begin{prop} In the above notation the following holds:

a) The groupoid $\Def _{\cR}(P)$ is connected.

b) The groupoid $\Def _{\cR}(I)$ is not connected.

c) The groupoid $\coDef _{\cR}(I)$ is connected.

d) The groupoid $\coDef _{\cR}(P)$ is not connected.
\end{prop}

\begin{proof} Let $(S,\id)\in \Def _{\cR}^{\h}(I)$. Then $S=I\otimes
_k\cR$ as a graded $(A\otimes \cR)^{op}$-module and the differential
in $S$ is equal to "multiplication by $\lambda (x\otimes
\epsilon)$" for some $\lambda \in k$. We denote this differential
$d_{\lambda}$ and the deformation $S$ by $S _{\lambda}$. By Lemma
11.7 each $(S_{\lambda}, \id)$ is also an object in the groupoid
$\Def _{\cR}(I)$. Notice that for $\lambda \neq 0$ we have
$H(S_{\lambda})=k$ and if $\lambda =0$ then
$H(S_{\lambda})=A\otimes \cR$. This shows for example that
$(S_1,\id)$ and $(S_0,\id)$ are non-isomorphic objects in $\Def
_{\cR}(I)$ and proves b).

The proof of d) is similar using Lemma 11.8.

Let us prove a). By Proposition 11.2, 1) the groupoid $\Def
_{\cR}(P)$ is equivalent to the full subcategory of $\Def
^{\h}_{\cR}(P)$ consisting of objects $(S,\id)$ such that $S\in
\cP(\cA _{\cR}^{op})$ or, equivalently, $i^*S=\bL i^*S$. As in the
proof of b) above we have $S=P\otimes \cR$ as a graded $(A\otimes
\cR)^{op}$-module and the differential in $S$ is equal to
"multiplication by $\lambda (x\otimes \epsilon)$" for some
$\lambda \in k$. Again we denote the corresponding $S$ by
$S_{\lambda}$. It is clear that the trivial homotopy deformation
$S_0$ is h-projective in $\cA ^{op}_{\cR}\text{-mod}$, hence it is
also an object in $\Def _{\cR}(P)$. It remains to prove that for
$\lambda \neq 0$ the DG $\cA _{\cR}^{op}$-module $S_{\lambda}$ is not
h-projective. Since the DG functor $\pi _!$ preserves
h-projectives (Example 3.13) it suffices to show that $S_{\lambda
}$ considered as a DG $\cR$-module is not h-projective. We have
$$\pi _!S_{\lambda}=\bigoplus_{n\geq 0}\cR [-n]$$
with the differential $\lambda \epsilon :\cR [-n]\to \cR [-n-1]$.
Consider the DG $\cR$-module
$$N=\bigoplus_{n=-\infty}^{\infty}\cR [-n]$$
with the same differential $\lambda \epsilon :\cR [-n]\to \cR
[-n-1].$ Note that $N$ is acyclic (since $\lambda \neq 0$) and the
obvious embedding of DG $\cR$-modules $\pi
_!S_{\lambda}\hookrightarrow N$ is not homotopic to zero. Hence $\pi
_!S_{\lambda}$ is not h-projective. This proves a).

The proof of c) is similar using Proposition 11.2, 2) and the DG
functor $\pi_*$ from Example 3.13.
\end{proof}

\end{document}